\documentclass[a4paper,12pt]{article}

\usepackage{amssymb,amsmath}
\usepackage{xypic}
\xyoption{all}
\usepackage{graphicx}

\newtheorem{defn}{Definition}[section]
\newtheorem{proposition}[defn]{Proposition}
\newtheorem{corollary}[defn]{Corollary}
\newtheorem{rem}[defn]{Remark}
\newtheorem{exm}[defn]{Example}
\newtheorem{lemma}[defn]{Lemma}
\newtheorem{theorem}[defn]{Theorem}
\newtheorem{notat}[defn]{Notation}
\newtheorem{newpar}[defn]{}

\newtheorem{xdefn}{Definition.}
\newtheorem{xproposition}{Proposition.}
\newtheorem{xcorollary}{Corollary.}
\newtheorem{xrem}{Remark.}
\newtheorem{xexm}{Example.}
\newtheorem{xlemma}{Lemma.}
\newtheorem{xtheorem}{Theorem.}
\newtheorem{xnotat}{Notation.}
\newtheorem{xnewpar}{\it}
\newtheorem{xproof}{{\it Proof. }}
\newtheorem{xproofof}{{\it Proof}}

\newenvironment{definition}{\begin{defn}\em}{\end{defn}}

\newenvironment{example}{\begin{exm}\em}{\end{exm}}

\newenvironment{proof}{\begin{xproof}\em}{\end{xproof}}

\newenvironment{newparagraph*}[1]{\begin{xnewpar}\hspace*{-1.5mm}{#1}. \rm}{\end{xnewpar}}

\newenvironment{definition*}{\begin{xdefn}\em}{\end{xdefn}}
\newenvironment{remark*}{\begin{xrem}\em}{\end{xrem}}
\newenvironment{example*}{\begin{xexm}\em}{\end{xexm}}
\newenvironment{notation*}{\begin{xnotat}\em}{\end{xnotat}}
\newenvironment{proposition*}{\begin{xproposition}}{\end{xproposition}}
\newenvironment{corollary*}{\begin{xcorollary}}{\end{xcorollary}}
\newenvironment{lemma*}{\begin{xlemma}}{\end{xlemma}}
\newenvironment{theorem*}{\begin{xtheorem}}{\end{xtheorem}}

\def\qed{\hspace{0.3cm}{\rule{1ex}{2ex}}}
\newcommand\V{\bigvee}
\newcommand\Max{\operatorname{Max}}
\newcommand\ie{i.e.}
\newcommand\eg{e.g.}

\newcommand\pwset[1]{\wp(#1)}
\newcommand\st{\mid}

\newcommand\cf{\textrm{cf.}}

\newcommand\opens{\operatorname{\mathcal{O}}}
\newcommand\spectrum{\operatorname{\Sigma}}
\newcommand\topology{\operatorname{\Omega}}
\newcommand\groupoid{\operatorname{\mathcal{G}}}
\newcommand\spp{\varsigma}
\newcommand\downsegment{{\downarrow}}

\newcommand\Frm{\textit{Frm}}
\newcommand\Loc{\textit{Loc}}
\newcommand\Top{\textit{Top}}
\newcommand\opp[1]{{#1}^{\textrm{op}}}
\newcommand\ootimes[1]{\mathop{\otimes}\limits_{#1}}
\newcommand\ptimes[1]{\mathop{\times}\limits_{#1}}
\newcommand\rs{\mathrm{R}}
\newcommand\ls{\mathrm{L}}
\newcommand\ident{\mathrm{id}}
\newcommand\ipi{\mathcal I}
\newcommand\lc{\operatorname{\mathcal L}}

\newcommand\CC{\mathbb{C}}
\newcommand\SL{\mathit{SL}}
\newcommand\lcc{\operatorname{{\mathcal L}^{\vee}}}
\newcommand\ssq{\mathit{StabQu}}
\newcommand\iqcat{\mathit{InvQu}}
\newcommand\iqf{\mathit{InvQuF}}
\newcommand\isgcat{\mathit{InvSGrp}}
\newcommand\apgcat{\mathit{ACPGrp}}
\newcommand\imcat{\mathit{InvMon}}
\newcommand\uiq{\mathit{Qu}}

\begin{document}

\title{\vspace*{-1cm}\'{E}tale groupoids and their quantales\thanks{Research supported in part by FEDER
and FCT/POCTI/POSI through grant POCTI/MAT/55958/2004 and the Center for Mathematical Analysis, Geometry, and Dynamical Systems, and by CRUP and British Council through the Treaty of Windsor grant No.\ B-22/04.}}
\author{Pedro Resende
\vspace*{2mm}\\ \small\it Departamento de Matem{\'a}tica, Instituto
Superior T{\'e}cnico, \vspace*{-2mm}\\ \small\it Av. Rovisco Pais 1,
1049-001 Lisboa, Portugal}

\date{~}

\maketitle
\vspace*{-1cm}
\begin{abstract}
We establish a close and previously unknown relation between quantales and groupoids, in terms of which the notion of \'{e}tale groupoid is subsumed in a natural way by that of quantale. In particular, to each \'{e}tale groupoid, either localic or topological, there is associated a unital involutive quantale. We obtain a bijective correspondence between localic \'{e}tale groupoids and their quantales, which are given a rather simple characterization and are here called inverse quantal frames. We show that the category of inverse quantal frames is equivalent to the category of complete and infinitely distributive inverse monoids, and as a consequence we obtain a correspondence between these and localic \'{e}tale groupoids that generalizes more classical results concerning inverse semigroups and topological \'{e}tale groupoids. This generalization is entirely algebraic and it is valid in an arbitrary topos.
As a consequence of these results we see that a localic groupoid is \'{e}tale if and only if its sublocale of units is open and its multiplication map is semiopen, and an analogue of this holds for topological groupoids. In practice we are provided with new tools for constructing localic and topological \'{e}tale groupoids, as well as inverse semigroups, for instance via presentations of quantales by generators and relations. The characterization of inverse quantal frames is to a large extent based on a new quantale operation, here called a support, whose properties are thoroughly investigated, and which may be of independent interest.
\vspace{0.2cm}
\vspace{0.2cm}
\end{abstract}

\maketitle

\newpage
\vspace*{1cm}
\tableofcontents
\newpage

\section{Introduction}\label{sec:introduction}

In this paper we study quantales, inverse semigroups, and groupoids. There are close and well known connections between the latter two concepts, whose importance in algebra, geometry, and analysis is by now firmly established (see, \eg, \cite{CannasWeinstein,Connes,Higgins,Lawson,Mackenzie2,Mackenzie,MoerdijkMrcun,Paterson,Renault}), and in this paper we intend to provide evidence showing that also the theory of quantales may have a natural role to play in this context, indeed providing useful new insights and techniques.

The historical motivation for this is clear. Quantales (short for ``quantum locales''~\cite{Mu86}) are algebraic structures that can be regarded as point-free noncommutative spaces (see, \eg, the surveys \cite{Mu02,PaRo00}), and it is therefore natural to investigate the extent to which they relate to basic notions in noncommutative geometry such as groupoids and operator algebras.

Several papers have been devoted to the latter, for instance proving early on that right sided quantales classify postliminary C*-algebras but not general C*-algebras~\cite{BoRoBo}. More recently, advances in the representation theory of quantales~\cite{MuPe01} have led to an improvement of this, namely enabling one to show that the quantale valued functor $\Max$ (a noncommutative analogue of the maximal spectrum of commutative C*-algebras~\cite{Mu89,MuPe02}) is a complete invariant of unital C*-algebras~\cite{KrRe}. Although this is a positive (and to some extent surprising) result, it is still unsatisfactory because $\Max$ does not have the functorial properties that would be appropriate for making it a noncommutative generalization of Gelfand duality~\cite{KrPeReRo}, the quest for which is one of the guiding ideas behind the theory.

This suggests that more specific examples should be looked at, such as C*-algebras obtained from locally compact groupoids~\cite{Paterson,Renault}, or just groupoids per se. Following this motivation, the first direct connection between a quantale and a (principal) groupoid, again via representation theory, appeared in a paper where the space of Penrose tilings has been modeled as a quantale~\cite{MuRe}.
The present paper is inspired by this but it is very different because (besides being concerned with more than the single example of Penrose tilings) it does not use representations, it also addresses inverse semigroups, and it is, from the point of view of the relation between quantales and groupoids, much broader in scope because it applies to arbitrary \'{e}tale groupoids and, to some extent, also to more general classes such as that of open groupoids.

An important consequence of our results is an equivalence (which extends more ``classical'' results concerning inverse semigroups and groupoids) between the following three concepts: localic \'{e}tale groupoids; complete and infinitely distributive inverse semigroups; and quantales of a kind known in this paper as \emph{inverse quantal frames}.

Sections \ref{sec:introgrpdqnts} and \ref{sec:mainresults} provide an overview of the main ideas and results of this paper, whose consequences and significance are discussed in section \ref{sec:significance}.

We shall assume, until the end of section 1, that the reader is familiar with the concepts of sup-lattice, frame, locale, quantale, inverse semigroup, and groupoid, hence giving little or no explanation for these terms, but in section 2 some background is provided.

\subsection{Preliminary ideas}\label{sec:introgrpdqnts}

Just as a $\CC$-algebra is a semigroup in the monoidal category of complex linear spaces, so a quantale is a semigroup in the category of sup-lattices, or, in other words, a sup-lattice equipped with an associative multiplication that distributes over arbitrary joins in both variables:
\begin{eqnarray*}
a\left({\V_i b_i}\right)&=& \V ab_i\\
\left({\V_i a_i}\right) b&=&\V a_i b\;.
\end{eqnarray*}
The groupoid quantale associated to a discrete groupoid $G$ is then, in analogy with groupoid algebras, the free sup-lattice on $G$ equipped with a quantale multiplication which is defined to be the ``bilinear'' extension of the multiplication of $G$, where the product of any noncomposable pair of arrows of $G$ is defined to be $0$.
Concretely, this is the powerset $\pwset G$ with multiplication calculated pointwise,
\[XY=\{xy\st x\in X,\ y\in Y,\ \textrm{$x$ and $y$ compose}\}\;,\]
or, equivalently, the set $2^G$ of boolean valued functions on $G$ with convolution multiplication:
\[(f*g)(x)=\V\{f(y)\wedge g(z)\st x=yz\}\;.\]
There is additional structure on these quantales, namely an involution given in $\pwset G$ by pointwise inversion,
\[X^*=\{x^{-1}\st x\in X\}\;,\]
or, in $2^G$, by the formula
\[f^*(x)=f(x^{-1})\;,\]
and a multiplicative unit consisting, in $\pwset G$, of the set $G_0$ of units of $G$, or, in $2^G$, of its characteristic function.

Such quantales are said to be unital and involutive, and examples of them are well known, namely the quantale $\pwset{X\times X}$ of binary relations on a set $X$~\cite{MuPe92}, or, for a group $G$, the quantale $\pwset G$, which is the image of $G$ under a left adjoint functor from involutive monoids to unital involutive quantales whose right adjoint is the forgetful functor, and which can equally be seen to be the image of $G$ under a functor from groups to unital involutive quantales whose right adjoint sends a unital involutive quantale $Q$ to its group of unitary elements $\mathcal U(Q)=\{u\in Q\st uu^*=u^*u=e\}$.

Let again $G$ be a groupoid. In addition to being a unital involutive quantale, $\pwset G$ is also a locale  whose points can be identified with the arrows of $G$, and the multiplication on $G$ can be reconstructed entirely from the multiplication in $\pwset G$: identifying the points of $\pwset{X\times X}$ with the atoms, two arrows $x$ and $y$ compose if and only if $\{x\}\{y\}\neq\emptyset$, in which case we have $\{xy\}=\{x\}\{y\}$. Hence, $G$ is recovered from $\pwset G$ up to an isomorphism.

The purpose of this paper is to establish a similar correspondence for more general groupoids, either topological or localic. For instance, for topological groupoids whose topologies are closed under pointwise multiplication, the topologies themselves become involutive quantales. 
This includes a very wide range of examples that arise in practice, namely all the topological groups, \'{e}tale groupoids, Lie groupoids, locally compact groupoids (in the sense of \cite{Paterson}), etc.
Similarly, any localic groupoid
whose multiplication map is at least semiopen gives rise to a quantale. For instance, any Grothendieck topos arises from such a localic groupoid \cite{JT}. We shall see that under reasonable assumptions the whole groupoid structure can be recovered from its ``quantal topology'', just as in the discrete case, in addition obtaining a precise characterization of those quantales that arise from \'{e}tale groupoids.

In the course of doing this we shall see that these quantales are intimately related to inverse semigroups; indeed, our results extend more ``traditional'' dualities between inverse semigroups and \'{e}tale groupoids. The quantales that mediate the extension are obtained from the inverse semigroups by a left adjoint functor, hence providing us with a precise definition of a sense in which the associated groupoids are ``universal''.

\subsection{Overview}\label{sec:mainresults}

Consider again the unital involutive quantale $\pwset G$ associated to a discrete groupoid $G$. The direct image of the domain map $d:G\to G$ is an operation \[\spp:\pwset G\to\pwset G\] that preserves unions and clearly satisfies the following properties, for all $U\subseteq G$:
\begin{eqnarray*}
\spp U&\subseteq& G_0\\
\spp U&\subseteq& UU^*\\
U&\subseteq&\spp UU\\
\spp(UG)&\subseteq&\spp U\;.
\end{eqnarray*}

In section 3 we shall thoroughly study unital involutive quantales equipped with an operation $\spp$ satisfying the first three properties, to which we call a \emph{support}. The existence of a support in a unital involutive quantale $Q$ implies, among other things, that the principal order ideal
\[\downsegment e=\{a\in Q\st a\le e\}\]
generated by the multiplicative unit of $Q$ is a locale that coincides with $\spp Q=\{\spp a\st a\in Q\}$ (lemma \ref{prop:supps}). We shall see that a support satisfying also the fourth axiom (such a support is said to be \emph{stable}) is unique if it exists, being given by any of the following formulas, among others (lemma \ref{lemma:defstab} and theorem \ref{thm:fullsubcat}):
\begin{eqnarray*}
\spp a &=& a1\wedge e\\
\spp a &=& aa^*\wedge e\\
\spp a &=& \bigwedge\{b\in \downsegment e\st a\le ba\}\;.
\end{eqnarray*}
We shall refer to any unital involutive quantale equipped with a (stable) support as \emph{(stably) supported}, and
we shall see that the category of stably supported quantales is full and reflective in the category of unital involutive quantales (theorem \ref{thm:fullreflective}). In particular, it follows that having a stable support is a property, rather than extra structure.

Still in section 3 it is seen (theorem \ref{thm:ipiQ}) that associated to each supported quantale $Q$ there is an inverse monoid $\ipi(Q)$ consisting of the \emph{partial units} of $Q$, by which we mean the elements $a\in Q$ satisfying $aa^*\le e$ and $a^*a\le e$. For instance, for the quantale $\pwset{X\times X}$ these are the partial bijections on $X$, and thus the semigroup associated to the quantale is the symmetric inverse monoid of $X$. This construction is functorial, and it has a converse: a functor from inverse monoids to stably supported quantales, which to each inverse monoid $S$ associates, in a universal way (the functor is left adjoint to the previous one), a quantale $\lc(S)$ that concretely consists of the set of downwards closed sets of $S$ (theorem \ref{thm:LS}). More, $\lc(S)$ is also an example of what we call an \emph{inverse quantale}, \ie, a supported quantale for which each element is a join of partial units. This concept plays a fundamental role in this paper. The inverse monoid $\ipi(Q)$ of an inverse quantale is an example of a complete and infinitely distributive inverse monoid,
and another adjunction is obtained relating these monoids and inverse quantales, where now the left adjoint assigns to such a monoid $S$ the set $\lcc(S)$ of downwards closed sets of $S$ that are also closed under the formation of joins of $S$ (corollary \ref{cor:adjunction}).

Section 4 investigates the properties of \emph{stable quantal frames}, which by definition are the stably supported quantales for which in addition the locale distributivity law
\[a\wedge\V_i b_i=\V_i a\wedge b_i\]
holds --- equivalently, the unital involutive quantales that satisfy locale distributivity and, in addition, the following two laws:
\begin{eqnarray*}
a1\wedge e&\le&aa^*\\
(a1\wedge e)a&\ge&a\;.
\end{eqnarray*}
This study culminates in a complete algebraic characterization (theorem \ref{inversiontheorem}) of those quantales for which there is an associated groupoid, which is then necessarily \emph{\'{e}tale}; that is, the domain map is a local homeomorphism. An important aspect of the characterization is given by a purely algebraic fact about quantales (lemma \ref{inversionlemma}), which states that a stable quantal frame is an inverse quantale (for stable quantal frames this just means that the maximum $1$ is a join of partial units) if and only if it satisfies the following laws:
\begin{eqnarray*}
(a\wedge e)1&=&\V_{xy^*\le a} x\wedge y\;,\\
1(a\wedge e)&=&\V_{x^*y\le a} x\wedge y\;.
\end{eqnarray*}
As we shall see, these laws relate directly to the inversion law of groupoids (lemma \ref{lemma:groupoidiffinversion}).

The following theorem summarizes the essential aspects of the characterization obtained, as well as the relation to groupoids, where $Q\otimes_{\spp Q} Q$ is the tensor product of $Q$ with itself over $\spp Q$ with respect to the $\spp Q$-$\spp Q$-bimodule structure of $Q$ which is defined by multiplication of elements of $\spp Q$ on the left and on the right, respectively, and where $\mu:Q\otimes_{\spp Q}Q\to Q$ is the quantale multiplication, which is well defined as a homomorphism on $Q\otimes_{\spp Q}Q$, rather than just $Q\otimes Q$, due to associativity.

\begin{theorem}
A stable quantal frame $Q$ coincides with the locale of arrows of a localic \'{e}tale groupoid $\groupoid(Q)$ if and only if the following equivalent conditions hold:
\begin{enumerate}
\item $Q$ is an inverse quantale;
\item $Q\cong\lcc(\ipi(Q))$.
\end{enumerate}
Furthermore, if these equivalent conditions hold then the right adjoint to the quantale multiplication,
\[\mu_*(a)=\V\{x\otimes y\st xy\le a\}\;,\]
is a join preserving map from $Q$ to $Q\otimes_{\spp Q}Q$,
the tensor product $Q\otimes_{\spp Q}Q$ coincides with the locale of composable pairs of arrows of $G$, and $\mu_*$ is the inverse image homomorphism of the groupoid multiplication. The support, seen as a map $\spp:Q\to\spp Q$, coincides with the direct image of the domain map of $\groupoid(Q)$.
\end{theorem}

A consequence of the results in section~4 is that we have been provided with a new way of constructing \'{e}tale groupoids, via quantales. For instance, it is possible to construct \'{e}tale groupoids algebraically by presenting their corresponding quantales by generators and relations. 
The natural setting for this, and that which the section addresses, is that of localic groupoids, and the results obtained are constructive in the sense of being valid in an arbitrary topos. Topological groupoids can be constructed in a canonical way from localic ones (see proposition \ref{prop:locgrpdtotopgrpd}), and because of this we are also provided with a means of constructing topological groupoids. Although the process of producing the point spectra of the locales involved raises constructivity issues (proving the existence of enough points usually requires the axiom of choice --- \cf\ \cite[\S II 3.4]{Johnstone}), at least this problem is left to the end.

Section 5 looks at a converse to section 4: constructing quantales from groupoids. 
To each localic \'{e}tale groupoid $G$ we associate an inverse quantal frame $\opens(G)$, and a (non-functorial) duality is obtained:

\begin{theorem}
Let $G$ be a localic \'{e}tale groupoid, and let $Q$ be an inverse quantal frame. Then we have isomorphisms
\begin{eqnarray*}
Q&\cong&\opens(\groupoid(Q))\\
G&\cong&\groupoid(\opens(G))\;.
\end{eqnarray*}
\end{theorem}

In addition, the following characterization of \'{e}tale groupoids, which depends to a large extent on lemma \ref{inversionlemma}, is found (corollary \ref{cor:quantalunitalisetale}):

\begin{theorem}\label{thm:intro2}
A localic groupoid $G$ is \'{e}tale if and only if its multiplication map is semiopen and the sublocale of units $G_0$ is open.
\end{theorem}

A similar characterization exists for topological groupoids (theorem \ref{thm:topetale}). In analogy with the localic case, it draws its inspiration from lemma \ref{inversionlemma}. However, section 5.3, where theorem \ref{thm:topetale} is proved, is entirely self contained, and can be read by anyone interested only in topological groupoids:

\begin{theorem}\label{thm:intro3}
Let $G$ be a topological groupoid. The following conditions are equivalent:
\begin{enumerate}
\item $G$ is \'{e}tale.
\item The unit space $G_0$ is open and the pointwise product of any two open sets of $G$ is an open set.
\item The unit space $G_0$ is open and the continuous domain map $d:G\to G_0$ is open.
\end{enumerate}
\end{theorem}

This sharpens previously known facts. For instance, based on \cite[Prop.\ 2.8]{Renault} we may
state:

\begin{corollary}
For a groupoid $G$ which is r-discrete in the sense of \cite[Def.\ 2.6]{Renault} (\ie, a locally compact Hausdorff groupoid whose unit space is open), the following conditions are equivalent:
\begin{enumerate}
\item $G$ admits a left Haar measure.
\item The pointwise product of any two open sets of $G$ is an open set.
\item The domain map of $G$ is open.
\end{enumerate}
\end{corollary}

Our results lead to the view of \'{e}tale groupoids as being essentially those groupoids whose topologies are unital quantales. Indeed, in the case of localic groupoids it is certainly appropriate to say that \'{e}tale groupoids \emph{are} quantales! The choice of point of view (groupoids or quantales) is not irrelevant, since the categories are not equivalent (\cf\ discussion at the end of section \ref{sec:locetfgrps}).

Another consequence of these results is that a precise bijective correspondence (up to isomorphisms) has been established between complete and infinitely distributive inverse semigroups, on one hand, and inverse quantal frames, on the other, namely given by the assignments $Q\mapsto\ipi(Q)$ and $S\mapsto\lcc(S)$. Indeed, we have:

\begin{theorem}\label{thm:intro4}
The category of complete and infinitely distributive inverse mo\-noids is equivalent to the category of inverse quantal frames.
\end{theorem}

This of course means that also a bijective correspondence (although not an equivalence of categories) has been established between these inverse semigroups and localic \'{e}tale groupoids. This provides a means of extending ``classical'' results concerning inverse semigroups and topological \'{e}tale groupoids to an entirely algebraic (and constructive) setting. The details of this extension will be dealt with in another paper, where in particular it will be seen that the locale points of an inverse quantal frame $\lcc(S)$ can be (no longer constructively) identified with ``germs'' of elements of $S$, matching well known constructions like the germ groupoid of a pseudogroup (see, \eg, \cite[Section 2.4]{Lawson} or \cite[Section 5.5]{MoerdijkMrcun}) or the universal groupoid of an inverse semigroup of \cite{Paterson}.

Another aspect of our results, as regards the relation between inverse semigroups and groupoids, is that since $\lcc(S)$ is defined by a universal property, it follows that the associated groupoid is in some sense a ``universal'' groupoid of $S$ (roughly, the topology of the groupoid, rather than the groupoid itself, is ``freely'' generated). This also suggests a new way of looking at the notion of universal groupoid of \cite{Paterson}, providing a perhaps clearer sense (as compared to that of \cite[Prop.\ 4.3.5]{Paterson}) in which such groupoids are universal.

\subsection{General comments}\label{sec:significance}

The quantales associated to \'{e}tale groupoids are also frames, and in principle it would have been possible to restrict to such quantales from the beginning, instead of providing a separate study of general supported quantales in section~3. However, the latter is justified, and indeed somewhat unavoidable, for two reasons. First, many of the good properties of stable quantal frames are entirely independent from frames, being rather consequences of the axioms of stably supported quantales. Second, working with general stably supported quantales is the natural thing to do when meets are not necessarily available, notably when dealing with arbitrary inverse semigroups as we often do in this paper. In this respect, an important role is played by the (non-immediate) fact that the category of stably supported quantales is full and reflective in the category of unital involutive quantales.

An immediate consequence of our results is, of course, that the language of quantales has been made available for studying \'{e}tale groupoids in a natural way, for instance leading to the characterizations described in theorems \ref{thm:intro2} and \ref{thm:intro3}, and to the relation with inverse semigroups expressed in theorem \ref{thm:intro4}. For instance, groupoids can be presented by generators and relations through quantales. More than that, \'{e}tale groupoids can be regarded as objects in the category of quantales, which for some applications may be a convenient category to work with. Furthermore, as already mentioned in section \ref{sec:mainresults}, this may also enable us to work constructively, by which is meant the possibility of interpreting definitions and theorems in an arbitrary topos. Far from being just a philosophical matter, a theorem proved constructively is in fact many theorems (an equivariant version, a version for bundles, etc.), and so constructivity is also a matter of practical relevance (\cf\  \cite{MuPe91}).

An important aspect of these results is that they suggest ways of generalizing the duality of inverse semigroups and \'{e}tale groupoids beyond \'{e}tale groupoids; the natural generalization of the notion of inverse semigroup is here provided by quantales, and it is certainly worth extending this kind of work so as to include open groupoids, groupoid representations, Morita equivalence, etc.

Our results also suggest new ways of assessing the extent to which quantales may be considered geometric objects, since groupoids are often regarded as generalized spaces in their own right --- this claim is supported by the role played by open localic groupoids in Joyal and Tierney's fundamental representation theorem for Grothendieck toposes~\cite{JT}, by the existence of homology and cohomology theories directly associated to \'{e}tale groupoids (see, \eg, \cite{CrMo}), by a wealth of applications of Lie groupoids in differential geometry \cite{CannasWeinstein,Connes,Mackenzie,MoerdijkMrcun}, etc.

In addition, specific topics worth looking at are: the \'{e}tale groupoid of Penrose tilings (in \cite{MuRe} the topology of the groupoid has been ignored); constructions of groupoids from C*-algebras, such as AF-groupoids from AF-algebras \cite{Renault}, or the dual groupoid of a C*-algebra of \cite{Idrissi} (see \cite{Renault:dualgroupoid}), along with their relations both to the kind of quantales studied in this paper and to the functor $\Max$; also applications in rather different areas, for instance in logic where supported quantales provide a convenient algebraic semantics for modal logic (see noninvolutive version in \cite{Catia}), furthermore providing a natural bridge to geometric structures such as, say, foliated manifolds via their holonomy groupoids.

Our results also provide examples of how a space equipped with continuous partial operations (such as a topological groupoid) can be regarded instead as a space whose additional structure, consisting of join preserving operations that are no longer partial, is entirely placed on the topology. This suggests an unusual way of looking at structured spaces in general. In particular, the points are secondary entities and the operations induced on them are in general ``fuzzy'' (\cf\ example \ref{fuzzyexample}).

\section{Preliminaries}

In this section we review some background and prove preliminary results, in addition establishing notation and terminology that is necessary for making the paper more self-contained. Readers who are familiar with the concepts addressed here may move directly to section \ref{sec:sq}, then referring back to this section as needed.

\subsection{Locales and quantales}\label{sec:locqu}

We begin by discussing sup-lattices, for which a general reference is \cite{JT}.

The category of \emph{sup-lattices} $\SL$ has as objects the complete lattices and as morphisms the maps that preserve arbitrary joins, which will be called \emph{homomorphisms}. The greatest element of a sup-lattice $L$ will be denoted by $1_L$, or just $1$, and the least element by $0_L$ or $0$.

Quotients of sup-lattices can be conveniently represented by closure operators: if $j:L\to L$ is a closure operator on a sup-lattice $L$ (\ie, a monotone map such that $\ident_L\le j$ and $j\circ j=j$) then the set of fixed points $L_j=\{x\in L \st j(x)=x\}$ is a sup-lattice under the same order as in $L$, and $j:L\to L_j$ is a surjective homomorphism. Any surjective homomorphism arises like this, up to isomorphism: if $h:L\to M$ is a surjective homomorphism, then  $M\cong L_j$ for the closure operator $j=h_*\circ h$ with $h_*$ the right adjoint of $h$.

The following well known fact will be used later:

\begin{proposition}\label{prop:supextension}
The forgetful functor from $\SL$ to the category of posets and monotone maps has a left adjoint that to each poset $P$ assigns the sup-lattice $\lc(P)$ of downwards closed subsets of $P$ (the $\lc$ stands for
``lower sets'')\footnote{We use the terminology ``downwards closed set'' or ``lower set'' instead of ``order ideal'' because order ideals are usually required to be nonempty.}. The unit $P\mapsto\lc(P)$ of the adjunction sends each $p\in P$ to the principal order ideal $\downsegment p=\{x\in P\st x\le p\}$.
\end{proposition}

\begin{proof}
In other words, we are stating that for each poset $P$ the sup-lattice $\lc(P)$ is a join completion of $P$ such that for each monotone map $f:P\to L$ to a sup-lattice $L$ there is a unique join preserving extension $\overline f:\lc(P)\to L$ of $f$ along $\downsegment(-):P\to\lc(P)$. This extension is defined, for each $U\in\lc(P)$, by $\overline f(U)=\V_{x\in U}f(x)$. \qed
\end{proof}

The category $\SL$ is monoidal, and the \emph{tensor product} $L\otimes M$ of two sup-lattices $L$ and $M$ is, similarly to abelian groups, the image of a universal bimorphism $L\times M\to L\otimes M$, where by a \emph{bimorphism} $f:L\times M\to N$ is meant a map that preserves joins in each variable separately.

Now we discuss quantales. For a general reference concerning basic algebraic facts about quantales see \cite{Rosenthal}.

A \emph{quantale} is a semigroup in $\SL$, \ie, a sup-lattice $Q$ equipped with an associative multiplication 
\[m:Q\otimes Q\to Q\;.\]
We write just $ab$ for the multiplication $m(a\otimes b)$ of two quantale elements.

Quantales are like rings whose underlying additive abelian groups have been replaced by sup-lattices, and there is a corresponding notion of module: by a \emph{left module} over a quantale $Q$ is meant a sup-lattice $M$ equipped with a left action of $Q$ in $\SL$, \ie, a sup-lattice homomorphism
\[\alpha:Q\otimes M\to M\]
that is associative with respect to $m$. We write $ax$ for $\alpha(a\otimes x)$.

A quantale $Q$ is \emph{unital} if it contains an element $e\in Q$ which is a multiplicative unit, \ie, for all $a\in Q$ we have
\[ea=ae=a\;,\]
and \emph{involutive} if it has an \emph{involution}, by which is meant a sup-lattice homomorphism $(-)^*:Q\to Q$ such that
\begin{eqnarray*}
a^{**}&=&a\\
(ab)^*&=&b^* a^*\;.
\end{eqnarray*}
An involution is said to be \emph{trivial} if it coincides with the identity map (in this case $Q$ is necessarily commutative).
A left module $M$ over a unital quantale $Q$ is \emph{unital} if $ex=x$ for all $x\in M$.

There are many examples of unital involutive quantales. We mention the following (only the first two are directly relevant to this paper):
\begin{enumerate}
\item The powerset of any discrete groupoid, as mentioned in section 1.1.
\item In particular, the quantale $\pwset{X\times X}$ of binary relations on a set $X$~\cite{MuPe92}.
\item The quantale of sup-lattice endomorphisms of an orthocomplemented sup-lattice $L$~\cite{MuPe92} (this is isomorphic to $\pwset{X\times X}$ if $L\cong\pwset X$) --- $L$ is both a right and a left module over the quantale.
\item More generally, the quantale of endomorphisms of a symmetric sup-lattice 2-form~\cite{Re04}.
\item The quantale $\operatorname{Sub}(R)$ of additive subgroups of a unital involutive ring $R$. The set of additive subgroups $\operatorname{Sub}(M)$ of any left $R$-module $M$ is a left module over
$\operatorname{Sub}(R)$.
\item The quantale $\Max A$ of norm-closed linear subspaces of a unital C*-algebra $A$. Any representation of $A$ on a Hilbert space $H$ yields a left module over $\Max A$ consisting of all the projections of $H$~\cite{KrRe,MuPe01,Re04}.
\end{enumerate}

A \emph{homomorphism} of quantales (resp.\ unital, involutive) is a sup-lattice homomorphism that is also a homomorphism of semigroups (resp.\ monoids, involutive semigroups), and a \emph{homomorphism} of left modules is a sup-lattice homomorphism which is equivariant with respect to the action. We shall denote by $\uiq$ the category of unital involutive quantales and their homomorphisms, giving no special name to the more general categories of unital quantales, involutive quantales, or just quantales.

An element $a\in Q$ of a quantale $Q$ is \emph{right-sided} if $a1\le a$, and \emph{left-sided} if $1a\le a$. The set of right-sided (resp.\ left-sided) elements is denoted by $\rs(Q)$ (resp.\ $\ls(Q)$). Then $\rs(Q)$ is a left module over $Q$, with action given by multiplication on the left. Similarly, there is a notion of \emph{right module} over $Q$, of which $\ls(Q)$ is an example.

For instance, a right-sided element of $\operatorname{Sub}(R)$ for a ring $R$ is precisely a right ideal of $R$, and a right-sided element of $\operatorname{Max} A$ for a C*-algebra $A$ is a norm-closed right ideal of $A$.

In \cite{JT} modules over commutative quantales were thoroughly studied in the context of a sup-lattice version of commutative algebra and descent theory. There is no standard reference for modules over arbitrary quantales. See, \eg, \cite{Re04} for more details, including modules over involutive quantales.

Quotients of quantales and modules can be represented by closure operators satisfying additional conditions. A \emph{nucleus} of quantales on a quantale $Q$ is a closure operator $j:Q\to Q$ such that
\[j(a)j(b)\le j(ab)\]
for all $a,b\in Q$ \cite{Rosenthal}. The set of fixed points $Q_j$ is then a quotient of quantales, where in $Q_j$ the multiplication is given by $(a,b)\mapsto j(ab)$, and every quotient of quantales arises like this up to isomorphism.

Similarly, if $M$ is a left $Q$-module, a \emph{nucleus} on $M$ is a closure operator $j:M\to M$ such that $aj(x)\le j(ax)$ for all $a\in Q$ and all $x\in M$.

This follows a general pattern. For instance, for involutive quantales the corresponding nuclei satisfy the condition $j(a)^*\le j(a^*)$, etc.

As another example of quantale we have the notion of \emph{frame}, by which is meant a sup-lattice satisfying the distributivity law
\[a\wedge\V_i b_i=\V_i a\wedge b_i\;,\]
of which the main example is the topology $\topology(X)$ of a topological space $X$. Accordingly, a \emph{basis} of a frame $L$ is a subset $B\subseteq L$ which is join-dense in the sense that every $x\in L$ is a join of elements of $B$:
\[x=\V\{b\in B\st b\le x\}\;.\]

Any frame is a unital involutive quantale with multiplication $\wedge$, $e=1$, and trivial involution $a^*=a$. A quantale is a frame if and only if it is unital with $e=1$ and idempotent~\cite{JT}.

The following simple result about frames will be useful later:

\begin{proposition}\label{prop:framebasis}
Let $h:L\to M$ be a frame homomorphism, and let $B\subseteq L$ be a basis of $L$ which furthermore is a downwards closed subset. Then $h$ is injective if and only if its restriction $h\vert_B:B\to M$ is.
\end{proposition}

\begin{proof}
Assume that $h\vert_B$ is injective. We shall prove that $h$ is injective (the converse is trivial).
Let $b\in B$ and $x\in L$ be arbitrary elements. Then,
\[\begin{array}{rcll}
h(b)\le h(x) &\iff& h(b)\wedge h(x)=h(b)\\
&\iff& h(b\wedge x)=h(b)\\
&\iff& b\wedge x=b & \textrm{(Because $b\wedge x\in B$.)}\\
&\iff& b\le x\;.
\end{array}
\]
Now let $x$ and $y$ be arbitrary elements of $L$. Then $y=\V Y$ for some $Y\subseteq B$, and we have
\[\begin{array}{rcll}
h(y)\le h(x) &\iff& h(\V Y)\le h(x)\\
&\iff& \V h(Y)\le h(x)\\
&\iff& \forall_{b\in Y}\ h(b)\le h(x)\\
&\iff& \forall_{b\in Y}\ b\le x\\
&\iff& y\le x\;.
\end{array}
\]
Hence, $h$ is an order embedding. \qed
\end{proof}

We denote by $\Frm$ the category of frames and their homomorphisms, and the opposite category $\opp\Frm$ by $\Loc$.
The identity functor on $\Frm$ can thus be seen as a contravariant functor
\[\opens:\opp{\Loc}\to\Frm\;,\]
to which we refer as the functor of \emph{opens}, in imitation of topological spaces --- but we shall keep the notations different, using $\opens$ for locales and $\topology$ for spaces. Following more or less standard terminology,
the objects of $\Loc$, which are the same as those of $\Frm$, will be called \emph{locales}, and the morphisms of $\Loc$ will be called (locale) \emph{maps}. Given a locale $L$, we shall usually write $\opens (L)$ instead of just $L$ when we wish to emphasize that we are thinking of $L$ as an object of $\Frm$. Similarly, we write $\opens(f)$, or $f^*$, for the frame homomorphism $\opens(M)\to\opens(L)$ corresponding to a map $f:L\to M$. 

If $L$ is a locale then by a \emph{quotient} of $\opens( L)$ is meant a \emph{nucleus} (of quantales) on $\opens (L)$, see~\cite{Johnstone}, or simply any surjective homomorphism defined on $\opens(L)$. This defines a \emph{sublocale} of $L$.

Again for a locale $L$, if $a\in \opens(L)$ then the map $(-)\wedge a:\opens(L)\to\downsegment a$ is a surjective frame homomorphism, and
$\downsegment a$ is said to be an \emph{open sublocale} of $L$.

The locale product $L\times M$ is defined by identifying $\opens(L\times M)$ with the coproduct of the frames $\opens(L)$ and $\opens(M)$, which can be conveniently described by their tensor product as sup-lattices,
\[\opens(L\times M)=\opens(L)\otimes\opens(M)\;,\]
with the binary meets calculated as follows:
\[(a\otimes b)\wedge(c\otimes d)=(a\wedge c)\otimes(b\wedge d)\;.\]

A locale map $f:L\to M$ is said to be \emph{semiopen}  if the frame homomorphism $f^*:\opens (M)\to\opens (L)$ preserves all the meets of $\opens(M)$; equivalently, if $f^*$ has a left adjoint $f_!$ --- the \emph{direct image} of $f$. (For adjoints between partial orders see \cite[Ch.\ IV]{MacLane} or \cite[Ch.\ I]{Johnstone}.)

A semiopen locale map $f:L\to M$ is \emph{open} if it satisfies the following condition for all $a\in\opens(L)$ and $b\in\opens(M)$,
\[f_!(a\wedge f^*(b))=f_!(a)\wedge b\;,\]
known as the \emph{Frobenius reciprocity condition} (see \cite[p.\ 521]{Elephant} or \cite[Ch.\ V]{JT}).
Equivalently~\cite[Ch.\ V]{JT}, the locale map $f$ is open if and only if it maps open sublocales to open sublocales, in the sense that the \emph{image} of each open sublocale $\downsegment a$, which is defined (up to isomorphism) to be the sublocale $X$ of $M$ determined by the regular epi-mono factorization
\[
\xymatrix{
\opens(M)\ar[rrrr]^{((-)\wedge a)\circ f^*}\ar@{->>}[drr]&&&&\downsegment a\\
&&\opens(X)\,,\ar@{>->}[urr]
}
\]
is open. (Then $X\cong\downsegment f_!(a)$.)
If $f$ is a continuous open map of topological spaces then $f^{-1}$, as a frame homomorphism, defines an open locale map.

An open locale map $f:L\to M$ is a \emph{local homeomorphism} if there is a cover $C$ of $L$ (a subset $C\subseteq L$ with $\V C=1$) such that for each $a\in C$ the frame homomorphism
\[((-)\wedge a)\circ f^*:\opens(M)\to\downsegment a\]
is surjective (this is the analogue for locales of a continuous open map of spaces whose restriction to each open set in a given cover is a subspace inclusion).

The assignment to each topological space $X$ of its locale of open sets $\topology(X)$ is the object part of a functor from $\Top$ to $\Loc$. This functor has a right adjoint
\[\spectrum:\Loc\to\Top\]
that to each locale $A$ assigns the space that best ``approximates'' it. Concretely, $\spectrum(A)$ is the set of points of $A$, \ie, the locale maps $p:1\to A$ where $1$ is the terminal object in $\Loc$, with a topology consisting of the open sets of the form $U_a=\{p\in\spectrum(A)\st p^*(a)=1\}$ for each $a\in\opens(A)$. The adjunction between $\Top$ and $\Loc$ restricts to an equivalence of categories between the \emph{sober spaces} and the \emph{spatial locales}.

For these and further facts about frames and locales see~\cite{Johnstone,JT}.

\subsection{Topological and localic groupoids}\label{sec:toplocgrpds}

In any category $C$ with pullbacks an \emph{internal groupoid} is a pair $(G_1,G_0)$, where $G_1$ is the \emph{object of arrows} and $G_0$ is the \emph{object of units}, equipped with two morphisms \[d,r:G_1\to G_0\;,\] called the \emph{domain} and \emph{range} morphisms, respectively, plus a \emph{multiplication} morphism \[m:G_1\ptimes{G_0} G_1\to G_1\;,\] where $G_1\times_{G_0} G_1$ is the pullback of $d$ and $r$ (the object of ``composable pairs of arrows''),
\[
\xymatrix{
G_1\ptimes{G_0} G_1\ar[rr]^-{\pi_1}\ar[d]_{\pi_2}&&G_1\ar[d]^r\\
G_1\ar[rr]_d&&G_0\;,
}
\]
a \emph{unit inclusion} morphism \[u:G_0\to G_1\;,\] and an \emph{inversion} morphism
\[i:G_1\to G_1\;,\]
all of which are required to satisfy the axioms (for associativity of multiplication, etc.) that make the data $(G_1,G_0,d,r,m,u)$ an internal category in $C$ (see \cite[Ch.\ XII]{MacLane}), and in addition the axioms for $i$, which consist of the equations
\begin{eqnarray*}
d\circ i &=& r  \\
r\circ i &=& d 
\end{eqnarray*}
and the inversion law expressed by the commutativity of the following diagram:
\[
\xymatrix{G_1\ar[d]_d\ar[rr]^-{\langle\ident,i\rangle}&&G_1\ptimes{G_0}G_1\ar[d]^m&&G_1\ar[ll]_-{\langle i,\ident\rangle}\ar[d]^r\\
G_0\ar[rr]_u&&G_1&&G_0.\ar[ll]^u}
\]

By a \emph{morphism} $f:G\to G'$ of internal groupoids in $C$ is meant a pair
\begin{eqnarray*}
f_1&:&G_1\to G'_1\\
f_0&:&G_0\to G'_0
\end{eqnarray*}
of morphisms of $C$ that commute in the natural way with the structure morphisms of the groupoids (\ie, $m'\circ f_1\times f_1=f_1\circ m$, $d'\circ f_1=f_0\circ d$, etc.\ --- \cf\ internal functors in \cite[Ch.\ XII]{MacLane}), a consequence of which is that in fact $f_0$ is uniquely determined by
$f_0=d'\circ f_1\circ u=r'\circ f_1\circ u$.

A straightforward application of the groupoid axioms shows that the inversion morphism $i$ is an involution:

\begin{proposition}\label{prop:invcat}
Let $G$ be an internal groupoid in a category $C$ with pullbacks. Then the inversion morphism $i$ satisfies the equations
\begin{eqnarray*}
i\circ i&=&\ident\\
i\circ m&=&m\circ\chi\;,
\end{eqnarray*}
where $\chi=\langle i\circ\pi_2,i\circ\pi_1\rangle$ is the unique morphism defined by the universal property of the pullback in the following diagram:
\[
\xymatrix{G_1\ptimes{G_0}G_1\ar@/_/[ddr]_{i\circ\pi_1}\ar@/^/[drrr]^{i\circ\pi_2}\ar@{.>}[dr]|\chi\\
&G_1\ptimes{G_0}G_1\ar[rr]^{\pi_1}\ar[d]^{\pi_2}&&G_1\ar[d]^r\\
&G_1\ar[rr]_d&&G_0}
\]
Furthermore $\chi$ is an isomorphism, and $\chi=\chi^{-1}$.
\end{proposition}

\begin{proof}
This is similar to the case of discrete groupoids, where the two conditions correspond of course to the equations $(x^{-1})^{-1} = x$ and $(xy)^{-1} = y^{-1}x^{-1}$.
Similarly to the discrete case, we apply the groupoid laws for $i$ together with associativity and the unit laws. For instance, for the first equation we have
\[\begin{array}{rcll}
i\circ i &=& m\circ\langle \ident,u\circ r\rangle\circ i\circ i &\textrm{(Unit law)}\\
&=& m\circ\langle i\circ i, u\circ r\circ i\circ i\rangle\\
&=& m\circ\langle i\circ i, u\circ r\rangle&(r\circ i=d,\ d\circ i=r)\\
&=& m\circ\langle i\circ i, m\circ\langle i,\ident\rangle\rangle&\textrm{(Inversion law)}\\
&=& m\circ\langle m\circ\langle i\circ i,i\rangle,\ident\rangle&\textrm{(Associativity)}\\
&=& m\circ\langle m\circ\langle i,\ident\rangle\circ i,\ident\rangle\\
&=& m\circ\langle u\circ r\circ i,\ident\rangle&\textrm{(Inversion law)}\\
&=& m\circ\langle u\circ d,\ident\rangle&(r\circ i=d)\\
&=& \ident\;. &\textrm{(Unit law)}
\end{array}\]
The second condition is equally straightforward, and the condition $\chi=\chi^{-1}$ follows from the pullback property, since this defines $\chi$ uniquely from the equations
$\pi_1\circ\chi=i\circ\pi_2$ and $\pi_2\circ\chi=i\circ\pi_1$,
which then lead to
\[\pi_1\circ\chi\circ\chi=i\circ\pi_2\circ\chi=i\circ i\circ\pi_1=\ident\circ \pi_1=\pi_1\circ\ident\]
and, similarly,
\[\pi_2\circ\chi\circ\chi=\pi_2\circ\ident\;,\]
hence showing that $\chi\circ\chi$ must coincide with $\ident$, again by the universal property of the pullback. \qed
\end{proof}

A \emph{localic groupoid (resp.\ category)} is an internal groupoid (resp.\ category) in $\Loc$, and a \emph{topological groupoid (resp.\ category)} is the same as an internal groupoid (resp.\ category) in $\Top$.

The axioms of groupoids include the equations $d\circ u=r\circ u=\ident_{G_0}$, and thus $G_0$ is always a subobject of $G_1$. In particular, for topological groupoids we shall often adopt the standard usage (see \cite{Paterson,Renault}) whereby a groupoid is thought of as a space $G$ equipped with suitable additional structure, which among other things implies the existence of the subspace $G_0$, defined to be the image $d(G)=r(G)$ of the domain and range maps, which are regarded as endomaps of $G$. Hence, expressions such as ``the topology $\topology(G)$ of $G$'', or ``the powerset $\pwset G$ of $G$'' have an obvious meaning.

A groupoid, either topological or localic, is said to be \emph{\'{e}tale} if the domain map $d$ is a local homeomorphism (this is equivalent to all the maps $d$, $r$, $m$, and $u$, being local homeomorphisms).

We remark that the notion of \emph{r-discrete} groupoid of \cite[Def.\ 2.6]{Renault} is closely related to \'{e}tale groupoids: it consists of a (Hausdorff, usually locally compact) groupoid $G$ whose unit space $G_0$ is open; the terminology is motivated by the fact that for an r-discrete groupoid $G$ the subspaces $r^{-1}(x)$, for all $x\in G_0$, are discrete. In practice, namely in the applications to operator algebras where a left Haar measure needs to be defined, r-discrete groupoids are necessarily \'{e}tale~\cite[Prop.\ 2.8]{Renault}. In \cite{Paterson} the definition of locally compact groupoid (no longer assumed Hausdorff) includes the existence of a left Haar measure, and the notions of r-discrete and \'{e}tale coincide.

\begin{proposition}\label{prop:locgrpdtotopgrpd}
The spectrum functor $\spectrum$ from locales to topological spaces induces a functor from the category of localic groupoids to the category of topological groupoids.
\end{proposition}

\begin{proof}
This a consequence of the fact that the spectrum functor $\spectrum$ preserves limits (because, as mentioned in section~\ref{sec:locqu}, it is right adjoint to the functor of opens $\topology$). In particular, for any localic groupoid $G$ the pullback
\[
\xymatrix{
G_1\ptimes{G_0} G_1\ar[rr]\ar[d]&&G_1\ar[d]^r\\
G_1\ar[rr]_d&&G_0\;,
}
\]
is mapped to a pullback of topological spaces
\[
\xymatrix{
\spectrum(G_1\ptimes{G_0} G_1)\ar[rr]\ar[d]&&\spectrum(G_1)\ar[d]^{\spectrum(r)}\\
\spectrum(G_1)\ar[rr]_{\spectrum(d)}&&\spectrum(G_0)\;,
}
\]
and thus $\spectrum(G_1\times_{G_0}G_1)$ is homeomorphic to $\spectrum(G_1)\times_{\spectrum(G_0)}\spectrum(G_1)$. Hence, from a localic groupoid
\[
\xymatrix{
G\ \ \ \ =\ \ \ \ G_1\ptimes{G_0}G_1\ar[r]^-m&G_1\ar@(ru,lu)[]_i\ar@<1.2ex>[rr]^r\ar@<-1.2ex>[rr]_d&&G_0\ar[ll]|u
}
\]
we obtain a topological groupoid
\[
\xymatrix{
\spectrum(G)\ \ \ \ =\ \ \ \ \spectrum(G_1)\ptimes{\spectrum(G_0)}\spectrum(G_1)\ar[rr]^-{\spectrum(m)}&&\spectrum(G_1)\ar@(ru,lu)[]_{\spectrum(i)}\ar@<1.2ex>[rr]^{\spectrum(r)}\ar@<-1.2ex>[rr]_{\spectrum(d)}&&\spectrum(G_0)\ar[ll]|{\spectrum(u)}
}
\]
because all the commutative diagrams that correspond to associativity of multiplication, etc., are similarly preserved. In an analogous way we see that any morphism of localic groupoids is mapped via $\spectrum$ to a morphism of topological groupoids, giving us a functor. \qed
\end{proof}

A functor in the opposite direction does not exist. In fact this is already the case for groups because the coproduct of spatial frames is not necessarily spatial~\cite[p.\ 61]{Johnstone}, and thus the multiplication map of a topological group $G$
\[m:G\times G\to G\]
has an inverse image
\[m^{-1}:\topology(G)\to\topology(G\times G)\]
whose codomain cannot in general be extended to $\topology(G)\otimes\topology(G)$.
In practice this is not a serious restriction because often the localic and the topological products coincide. For instance, if a topological space $X$ is locally compact we have $\topology(X\times Y)\cong\topology(X)\otimes\topology(Y)$ for any topological space $Y$~\cite[p.\ 61]{Johnstone}, and thus, for instance, locally compact groups yield localic groups.

Let us look more closely at the pullback locale $G_1\times_{G_0} G_1$ of a localic groupoid (in fact the following comments apply to any pullback of locales). As a frame, this is defined by the pushout
\[
\xymatrix{\opens(G_0)\ar[rr]^{d^*}\ar[d]_{r^*}&&\opens(G_1)\ar[d]^-{\pi_2^*}\\
\opens(G_1)\ar[rr]_-{\pi_1^*}&&\opens(G_1)\ootimes{\opens(G_0)}\opens(G_1)\;,}
\]
where $\opens(G_1)\otimes_{\opens(G_0)}\opens(G_1)$ is a quotient of the coproduct $\opens(G_1)\otimes\opens(G_1)$ in $\Frm$. As remarked in section \ref{sec:locqu}, this coproduct coincides with the tensor product of sup-lattices, with meet calculated as
\[(a\otimes b)\wedge (c\otimes d)=(a\wedge c)\otimes(b\wedge d)\;.\]
Then the quotient is defined by the condition
\[1\otimes d^*(a) = \pi_2^*(d^*(a))=\pi_1^*(r^*(a))=r^*(a)\otimes 1\]
which stabilizing under meets becomes, for all $b,c\in\opens(G_1)$ and $a\in\opens(G_0)$:
\[
b\otimes (d^*(a)\wedge c) = (b\wedge r^*(a))\otimes c\;.
\]

\subsection{Inverse semigroups}\label{sec:prelim.invsems}

We shall now give an overview of basic inverse semigroup theory, in the course of which we shall prove some simple results that will be needed in this paper. As general references for inverse semigroups we suggest \cite{Lawson} or \cite{Paterson}.

By an \emph{inverse} of an element $x$ of a semigroup is meant an element $y$ in the semigroup such that
\begin{eqnarray*}
xyx &=& x\\
yxy&=&y\;.
\end{eqnarray*}
An \emph{inverse semigroup} is a semigroup for which each element has a unique inverse. Equivalently, an inverse semigroup is a semigroup for which each element has an inverse (hence, a regular semigroup) and for which any two idempotents commute.

In an inverse semigroup the inverse operation defines an involution, and we shall always denote the inverse of an element $x$ by $x^{-1}$ or $x^*$. The set of idempotents of an inverse semigroup $S$ is denoted by $E(S)$.

An \emph{inverse monoid} is an inverse semigroup that has a multiplicative unit, which is usually denoted by $e$.

A semigroup homomorphism beween inverse semigroups automatically preserves inverses. The  category whose objects are the inverse semigroups and whose arrows are the semigroup homomorphisms will be denoted by $\isgcat$. The category whose objects are the inverse monoids and whose arrows are the monoid homomorphisms will be denoted by $\imcat$.

As examples we have:
\begin{enumerate}
\item The set $\ipi(X)$ of partial bijections on a set $X$ (\ie, bijections from a subset of $X$ onto another subset of $X$) with the multiplication $fg$ of partial bijections $f:U\to V$ and $g:U'\to V'$ being the usual composition $g\circ f$ defined on the set $f^{-1}(V\cap U')\subseteq U$, with codomain $g(V\cap U')\subseteq V'$ (this is an inverse monoid). The inverse of a partial bijection is the inverse of a function in the usual sense.
\item More generally, any set $S$ of partial homeomorphisms on a topological space, where $S$ is closed under the above multiplication and inverses, is an inverse semigroup. This is called a \emph{pseudogroup}.
\item Any subsemigroup of operators on a Hilbert space consisting entirely of partial isometries and closed under adjoints.
\item The set $\ipi(G)$ of ``G-sets'' (in the sense of~\cite[p.\ 10]{Renault}) of a discrete groupoid $G$, under pointwise multiplication and inverses, where by a \emph{G-set} is meant a set $U$ for which both the domain and range maps are injective when restricted to $U$. (We remark that this terminology is unfortunate because it colides with the standard usage of ``$G$-set'' for a set equipped with an action by a group $G$ --- see, \eg, \cite{MacLaneMoerdijk}.) This is an inverse monoid with unit $G_0$.
\end{enumerate}

The Wagner-Preston theorem asserts that every inverse semigroup is concretely representable as a pseudogroup.

The \emph{natural order} of an inverse semigroup $S$ is a partial order, defined as follows:
\[x\le y\iff x=fy\textrm{ for some }f\in E(S)\;.\]
The product of any two idempotents $f$ and $g$ is their meet, $fg=f\wedge g$, and
$S$ is an inverse monoid if and only if the set of idempotents has a join, in which case we have
$e=\V E(S)$.
In the case of a pseudogroup an idempotent $f$ is the identity map on an open set $U$, and thus the natural order becomes
\[x\le y\iff x=y\vert_U\textrm{ for some open set }U\subseteq\operatorname{dom}(y)\textrm{ such that }\ident_U\in S\;.\]
Hence, the natural order is just the restriction order on partial maps.

An important property of any pseudogroup of the form $P=\ipi(X)$ is that it is distributive in the sense that the multiplication distributes over all the joins that exist; that is, for all $s\in P$ and all $X\subseteq P$ such that $\V X$ exists in $P$ we have that both $\V_{x\in X} sx$ and $\V_{x\in X} xs$ exist in $P$ and
\[s(\V X)=\V_{x\in X} sx\ \ \textrm{ and }\ \ (\V X)s = \V_{x\in X} xs\;.\]
Accordingly, we adopt the following definition:

\begin{definition}
An inverse semigroup $S$ is said to be \emph{infinitely distributive} if, for all $s\in S$ and all subsets $X\subseteq S$ for which $\V X$ exists in $S$, the following conditions hold:
\begin{enumerate}
\item $\V_{x\in X} sx$ exists in $S$;
\item $\V_{x\in X} xs$ exists in $S$;
\item $s(\V X)=\V_{x\in X} sx$;
\item $(\V X)s = \V_{x\in X} xs$.
\end{enumerate}
\end{definition}

\begin{remark*}
(This paragraph may be skipped by readers not interested in matters of constructivity.)
Infinite distributivity is defined only for $X\neq\emptyset$ in \cite[p.\ 28]{Lawson}. This restriction is to a large extent irrelevant because if the least element $0=\V\emptyset$ exists (this is an idempotent) then for all $s\in S$ we necessarily have $0s=s0=0$ because $0\le 0s\Rightarrow 0=(00^{-1})0s=0s$ and $s0=(0s^{-1})^{-1}=0$; that is, $0$ is preserved by multiplication, and thus the definition of infinite distributivity does not seem to depend on whether we restrict $X$ to be nonempty or not. However, although this is certainly true in classical mathematics, it may no longer be true internally in a topos whose internal logic does not satisfy the law of excluded middle: stating an assertion separately for $X=\emptyset$ and $X\neq\emptyset$ is in general not equivalent to stating it for all $X$. Since we want all the definitions and proofs in this paper to carry over to an arbitrary topos we shall take care of getting rid of restrictions that lead to unnecessary ``case analyses'', which is why we have adopted the stronger notion of infinite distributivity, even though classically the two notions are of course equivalent.
\end{remark*}

A very important property of inverse semigroups is the following:

\begin{proposition}
Let $S$ be an inverse semigroup. The following conditions are equivalent:
\begin{enumerate}
\item $S$ is infinitely distributive;
\item $E(S)$ is infinitely distributive.
\end{enumerate}
\end{proposition}

\begin{proof}
This is proved in \cite{Lawson} for infinite distributivity with respect to joins of non-empty sets, but the proof applies equally to joins of any subsets. \qed
\end{proof}

An analogous and related property of inverse semigroups concerns distributivity of binary meets over joins:

\begin{proposition}[Resende \cite{Re:SF}]
Let $S$ be an infinitely distributive inverse semigroup, let $x\in S$, and let $(y_i)$ be a family of elements of $S$. Assume that the join $\V_i y_i$ exists, and that the meet $x\wedge\V_i y_i$ exists. Then,
for all $i$ the meet $x\wedge y_i$ exists, the join $\V_i(x\wedge y_i)$ exists, and we have
\[x\wedge\V_i y_i = \V_i(x\wedge y_i)\;.\]
\end{proposition}

Hence, the distributivity in $E(S)$ determines the distributivity in the whole of $S$, both with respect to multiplication and to binary meets (which in $E(S)$ are the same, of course).

We remark that the Wagner-Preston theorem gives us, for every inverse semigroup $S$, an injective homomorphism of inverse semigroups $r:S\to\ipi(X)$, where $X$ is a topological space. The representation $r$ is necessarily a monotone map, but it is important to notice that in general it does not preserve joins. In order to see this it suffices to consider a frame $L$ (this is an inverse semigroup with $E(L)=L$); then an injective homomorphism $r:L\to\ipi(X)$ that preserves arbitrary joins exists if and only if $L$ is spatial.

Finally, we recall that a pseudogroup $P$ of partial homeomorphisms on a topological space $X$ is said to be \emph{complete} if, for every partial homeomorphism $h$ on $X$ and every open cover $(U_i)$ of the domain of $h$, we have $h\in P$ whenever $h\vert_{U_i}\in P$ for all $i$. This can be equivalently formulated in terms of compatible elements of $P$, by which we mean any two elements $s,t\in P$ that coincide on the intersection of their domains. Noticing that $s$ and $t$ are compatible if and only if both $st^{-1}$ and $s^{-1}t$ are idempotents in $P$, we are led to the following definitions:

\begin{definition}
Let $S$ be an inverse semigroup. Two elements $s,t\in S$ are said to be \emph{compatible} if both $st^{-1}$ and $s^{-1}t$ are idempotents. A subset $X\subseteq S$ is \emph{compatible} if any two elements in $X$ are compatible. Then $S$ is said to be \emph{complete} if every compatible subset $X$ has a join $\V X$ in $S$ (hence, $S$ is necessarily a monoid with $e=\V E(S)$).
\end{definition}

We remark that we have defined completeness with respect to arbitrary compatible subsets (instead of just non-empty ones as in \cite{Lawson}). Hence, a complete inverse semigroup necessarily has a least element $0$.

As mentioned above, pseudogroups of the form $\ipi(X)$ are infinitely distributive, but indeed this is a property of any complete pseudogroup. Accordingly, in order to simplify our terminology we adopt the following (nonstandard) definitions:

\begin{definition}
By an \emph{abstract complete pseudogroup} will be meant a complete and infinitely distributive inverse semigroup (hence a monoid). The \emph{category of abstract complete pseudogroups} $\apgcat$ has the abstract complete pseudogroups as objects and the monoid homomorphisms that preserve the joins of all the compatible sets as arrows.
\end{definition}

We list here some properties of abstract complete pseudogroups that will be useful later on:

\begin{proposition}\label{prop:abspg}
Let $S$ and $T$ be abstract complete pseudogroups. Then,
\begin{enumerate}
\item $E(S)$ is a frame;
\item $S$ is meet-complete; that is, if $X\neq\emptyset$ then $\bigwedge X$ exists in $S$;
\item A monoid homomorphism $h:S\to T$ is in $\apgcat$ if and only if the restriction $h\vert_{E(S)}$ preserves joins (equivalently, $h\vert_{E(S)}$ is a homomorphism of frames);
\item The (obvious) forgetful functor from $\apgcat$ to $\imcat$ has a left adjoint.
\end{enumerate}
\end{proposition}

\begin{proof}
1. The first condition is obvious.

2. In order to prove the second condition, which as we shall see depends only on the fact that $E(S)$ is a frame,
define the following idempotent:
\[f=\V\{g\in E(S)\st gx=gy\textrm{ for all }x,y\in X\}\;.\]
Then, due to infinite distributivity, $fs=ft$ for all $s,t\in X$:
\begin{eqnarray*}
fs&=&\V\{gs\st gx=gy\textrm{ for all }x,y\in X\}\\
&=&\V\{gt\st gx=gy\textrm{ for all }x,y\in X\}\\
&=&ft\;.
\end{eqnarray*}
Let us denote by $z$ the (unique) value $fx$ for $x\in X$. This is a lower bound of the set $X$. Let $w$ be another lower bound. Then $w=ww^{-1}x$ for all $x\in X$, and thus $ww^{-1}\le f$. Hence, $w\le z$, which shows that $z=\bigwedge X$.

3. Let $h:S\to T$ be a homomorphism of monoids, and assume that $h\vert_{E(S)}$ preserves joins. We shall prove  that $h$ preserves all the joins that exist in $S$ (the converse implication is trivial).
Let $X\subseteq S$ be a subset for which the join $\V X$ exists. For all $x\in X$ we have $x=xx^{-1}\V X$, and thus
$\V X=\V_x (xx^{-1}\V X)$. Then, due to infinite distributivity,
\[\V X=\left({\V_x (xx^{-1})}\right)\V X\;,\]
and thus, using again infinite distributivity,
\[h\left({\V X}\right)=h\left({\left({\V_x (xx^{-1})}\right)\V X}\right)=h\left({\V_x (xx^{-1})}\right)h\left({\V X}\right)=\]
\[\left({\V_x h(xx^{-1})}\right) h\left({\V X}\right)=\V_x \left({h(x)h(x)^{-1} h\left({\V X}\right)}\right)=
\V_x h(x)\;.\]

4. The fourth property is a restriction (to monoid homomorphisms) of the adjunction of \cite[Section 1.4]{Lawson}, which takes place between $\isgcat$, on one hand, and the category of complete and infinitely distributive inverse semigroups with join-preserving semigroup homomorphisms, on the other (with an obvious minor adaptation due to the definition of completeness in \cite{Lawson}, which concerns only 
 joins of non-empty sets). \qed
\end{proof}

\section{Supported quantales}\label{sec:sq}

In this section we study the notion of supported quantale that was alluded to in section \ref{sec:mainresults}. First, in section \ref{sec:basicsupport}, we develop the basic notion of support, and then in section \ref{sec:stablesupports} we address stable supports. The basic properties of supports and stable supports have a clear intuitive meaning if one keeps in mind as an example the powerset $\pwset G$ of a discrete groupoid $G$, and they all follow from a very simple set of axioms. Then in sections \ref{sec:partialunits}, \ref{sec:inversequantales}, and \ref{sec:envelopes} we study the relations between supported quantales and inverse semigroups.

\subsection{Supports}\label{sec:basicsupport}

\begin{definition}
Let $Q$ be a unital involutive quantale. A \emph{support} on $Q$ is a sup-lattice endomorphism
\[{\spp}:Q\to Q\]
satisfying, for all $a\in Q$,
\begin{eqnarray}
\spp a&\le& e\label{def:s1}\\
\spp a&\le& a  a^*\label{def:s2}\\
a&\le& \spp a a\label{def:s3}
\end{eqnarray}
A \emph{supported quantale} is a unital involutive quantale equipped with a specified support.
\end{definition}

\begin{example}
\begin{enumerate}
\item The only support of a frame is the identity.
\item A support for the powerset $\pwset G$ of a discrete groupoid $G$ is obtained from the domain map of $G$:
\[\spp U=\{d(x)\st x\in U\}\;.\]
As we shall see in section \ref{sec:stablesupports}, this is the only possible support for such a quantale.
\item In particular, the support of the quantale $\pwset{X\times X}$ of binary relations on a set $X$ is given by
$\spp R=\{(x,x)\st (x,y)\in R\textrm{ for some }y\}$.
\end{enumerate}
\end{example}

Now we address general properties of supported quantales.

\begin{lemma}\label{prop:supps}
Let $Q$ be a supported quantale. The following conditions hold:
\begin{eqnarray}
\spp a &=& a, \textrm{ for all }a\le e     \label{prop:s1}\\
\spp\spp a &=&\spp a     \label{prop:s2}\\
a &=& \spp b a, \textrm{ if } \spp a\le \spp b\label{prop:s3a}\\
a &=& \spp a a     \label{prop:s3}\\
(\spp a)^* &=& \spp a     \label{prop:s7a}\\
a &=& a\spp(a^*)     \label{prop:s4}\\
\spp a=0 &\Leftrightarrow&a=0     \label{prop:s4a}\\
\spp a&\le&\spp(a a^*)     \label{prop:s5}\\
\spp a 1 &=& a 1     \label{prop:s6}\\
a 1&=&a a^* 1     \label{prop:s6a}\\
\spp a &=&\spp a\spp a     \label{prop:s7}\\
a &\le& a a^* a     \label{prop:s8}\\
\spp(a 1) b&=&a1\wedge b  \label{prop:s9}\\
\spp(a 1)&=&a1\wedge e\label{prop:s10}\\
\spp(a\wedge b)&\le&a b^*\label{prop:s10a}
\end{eqnarray}
Furthermore,
\begin{itemize}
\item $Q$ is a Gelfand quantale, in the sense of~\cite{MuPe92};
\item the subquantale $\downsegment e$ coincides with $\spp Q$ and it is a locale with $a b=a\wedge b$;
\item all the elements of $\spp Q$ are projections, and $\spp Q$ is a unital involutive subquantale with trivial involution;
\item the sup-lattice homomorphism $\spp Q\to\rs(Q)$ defined by $a\mapsto a 1$ is a retraction split by the map
$\rs(Q)\to\spp Q$ which is defined by $a\mapsto\spp a$;
\item the map $\spp:\rs(Q)\to\spp Q$ is an order embedding.
\end{itemize}
\end{lemma}

\begin{proof}
First we prove properties (\ref{prop:s1})--(\ref{prop:s10a}).
\begin{itemize}

\item[(\ref{prop:s1}):] From (\ref{def:s3}) and (\ref{def:s1}), if $a\le e$ we have $a\le\spp a a\le\spp a e=\spp a$, and from (\ref{def:s2}) and (\ref{def:s1}) we have $\spp a\le a a^*\le a e^*=a e= a$.

\item[(\ref{prop:s2}):] Immediate from the previous one because $\spp a\le e$.

\item[(\ref{prop:s3a},\ref{prop:s3}):] From (\ref{def:s3}) and (\ref{def:s1}): if $\spp a\le \spp b$ we have
$a\le\spp a a\le\spp b a\le ea=a$.

\item[(\ref{prop:s7a}):] We have
$\spp a=\spp \spp a\le (\spp a) (\spp a)^*\le e(\spp a)^*$, and thus $\spp a\le (\spp a)^*$, \ie, $\spp a=(\spp a)^*$.

\item[(\ref{prop:s4}):] From (\ref{prop:s3}) and (\ref{prop:s7a}) we have $a=a^{**}= (\spp(a^*) a^*)^*=a\spp(a^*)^*=a\spp(a^*)$.

\item[(\ref{prop:s4a}):] If $\spp a=0$ then $a=0$ because $a\le\spp a a$. The converse, \ie, $\spp 0=0$, is trivial because $\spp$ preserves joins (but we remark that the axiom $\spp a\le a a^*$ would also imply $\spp 0=0$ for more general maps $\spp$).

\item[(\ref{prop:s5}):] Follows from (\ref{def:s2}) and (\ref{prop:s2}).

\item[(\ref{prop:s6},\ref{prop:s6a}):] Follows from (\ref{def:s2}) and (\ref{def:s3}):
$\spp a 1\le a a^* 1\le a 1\le \spp a a 1\le \spp a 1$.

\item[(\ref{prop:s7}):] Follows from (\ref{def:s3}) and (\ref{prop:s2}):
$\spp a\le\spp\spp a\spp a=\spp a \spp a\le\spp a$.

\item[(\ref{prop:s8}):] Follows from (\ref{def:s3}) and (\ref{def:s2}): $a\le\spp a a\le a a^* a$.

\item[(\ref{prop:s9}):] From (\ref{prop:s6}) we have
$\spp(a1)b\le\spp(a1)1=a11=a1$. Since
$\spp(a1)b\le eb=b$, we obtain the inequality
\[\spp(a1)b\le a1\wedge b\;.\]
The converse inequality follows from (\ref{def:s3}):
\[a1\wedge b\le\spp(a1\wedge b)(a1\wedge b)\le\spp(a1)b\;.\]

\item[(\ref{prop:s10}):] Follows from the previous one with $b=e$.

\item[(\ref{prop:s10a}):] Follows from (\ref{def:s2}):
$\spp(a\wedge b)\le(a\wedge b)(a\wedge b)^*\le ab^*$.

\end{itemize}
A Gelfand quantale is one for which $a=a a^* a$ for all right-sided elements $a$. Hence, $Q$ is Gelfand because, if $a$ is right-sided, we have $a a^* a\le a$, which together with (\ref{prop:s8}) implies the Gelfand condition.

The downsegment $\downsegment e$ coincides with $\spp Q$ due to (\ref{prop:s1}). It
 is of course a unital subquantale, and it is idempotent due to (\ref{prop:s7}). Therefore it is an idempotent quantale whose unit is the top, in other words a locale with $a b=a\wedge b$.

We have already seen that the elements $a\le e$ are idempotent and, by (\ref{prop:s7a}), self-adjoint, \ie, projections. Hence, the locale $\spp Q$, with the trivial involution, is an involutive subquantale of $Q$.

Now we verify that the support splits the map $(-) 1:\spp Q\to\rs(Q)$. Let $a\in\rs(Q)$. Then,
by (\ref{prop:s6}), $\spp a 1=a 1=a$.

It follows that $\spp:\rs(Q)\to\spp Q$ is an order embedding because it is a section. \qed
\end{proof}

\subsection{Stable supports}\label{sec:stablesupports}

\begin{lemma}\label{lemma:defstab}
Let $Q$ be a supported quantale. The following conditions are equivalent:
\begin{enumerate}
\item\label{prop:s11} for all $a,b\in Q$, $\spp(a b) =\spp(a\spp b)$;
\item\label{prop:s12} for all $a,b\in Q$, $\spp(a b) \le \spp a$;
\item\label{prop:s14} for all $a\in Q$, $\spp(a 1)=\spp a$;
\item\label{prop:s15} for all $a\in Q$, $a1\wedge e=\spp a$;
\item\label{prop:s15a} for all $a,b\in Q$, $a1\wedge b=\spp a b$;
\item\label{prop:ss2} the map $(-) 1:\spp Q\to\rs(Q)$ is an order isomorphism whose inverse is $\spp$ restricted to $\rs(Q)$ (in particular, $\rs(Q)$ is a frame with $a1\wedge b1=\spp a b1$);
\item\label{prop:ss3} for all $a,b\in Q$, $a\le\spp b a$ if and only if $\spp a\le\spp b$;
\item\label{prop:ss4} for all $a\in Q$, $\spp a=\bigwedge\{b\in\spp Q\st a\le b a\}$;
\item\label{prop:ss5} for all $a\in Q$, $\spp a\le b 1$ if and only if $\spp a\le\spp b$;
\item\label{prop:ss6} for all $a\in Q$, $\spp a=\bigvee\{b\in\spp Q\st b\le a 1\}$;
\item\label{prop:ss7} the sup-lattice bimorphism $Q\times\spp Q\to\spp Q$ defined by $(a,f)\mapsto\spp(a f)$ makes $\spp Q$ a left $Q$-module (in this case the isomorphism $\spp Q\cong\rs(Q)$ of \ref{prop:ss2} is also a module isomorphism).
\end{enumerate}
\end{lemma}

\begin{proof}
First we show that the first five conditions are equivalent.
First, assuming \ref{prop:s11}, we have $\spp(a b)=\spp(a\spp b)\le\spp(a e)=\spp a$, which proves \ref{prop:s12}. Conversely, if \ref{prop:s12} holds then
\[\spp(a b)\le\spp(a\spp b b)\le\spp(a\spp b)\le\spp(a b b^*)\le\spp(a b)\;,\]
and thus \ref{prop:s11} holds. It is obvious that \ref{prop:s14} is equivalent to \ref{prop:s12}, and now let us show that \ref{prop:s14}, \ref{prop:s15}, and \ref{prop:s15a} are equivalent. First, (\ref{prop:s9}) tells us that $\spp(a1)b=a1\wedge b$, and thus if \ref{prop:s14} holds we obtain $\spp a b=a1\wedge b$. Hence, \ref{prop:s14} implies \ref{prop:s15a}, which trivially implies \ref{prop:s15}. Finally, if \ref{prop:s15} holds we have $\spp (a1)=a11\wedge e=a1\wedge e=\spp a$, and thus \ref{prop:s14} holds.

Now we deal with the remaining conditions.

(\ref{prop:s11}$\Rightarrow$\ref{prop:ss2}) Assume \ref{prop:s11}, and let $a\in Q$. Then
$\spp(\spp a 1)=\spp(\spp a \spp 1)=\spp(\spp a e)=\spp\spp a=\spp a$. This shows that the map $(-) 1:\spp Q\to\rs(Q)$ is an order isomorphism with $\spp$ as its inverse, because we have already seen that it is a retraction split by $\spp$.

(\ref{prop:ss2}$\Rightarrow$\ref{prop:s14}) Let $a\in Q$. We have $a 1=\spp a 1$ for any support, and thus
assuming \ref{prop:ss2} we have $\spp(a 1)=\spp(\spp a 1)=\spp a$.

(\ref{prop:s12}$\Rightarrow$\ref{prop:ss3}) Assume $a\le\spp b a$ and that \ref{prop:s12} holds. Then $\spp a\le\spp(\spp b a)\le\spp\spp b=\spp b$. The converse is the condition that $a\le\spp b a$ follows from $\spp a\le\spp b$, which coincides with (\ref{prop:s3a}).

(\ref{prop:ss3}$\Rightarrow$\ref{prop:s12}) From the condition $a\le\spp a a$ we obtain, multiplying by $b$ on the right,
$(a b)\le\spp a (a b)$, and thus assuming \ref{prop:ss3} we obtain $\spp(a b)\le\spp a$.

(\ref{prop:ss3}$\Rightarrow$\ref{prop:ss4}) If \ref{prop:ss3} holds then in the trivial equation $\spp a=\bigwedge\{b\in\spp Q\st\spp a\le b\}$ we can substitute $a\le b a$ for the equivalent condition $\spp a\le b$, thus obtaining $\spp a=\bigwedge\{b\in\spp Q\st a\le b a\}$.

(\ref{prop:ss4}$\Rightarrow$\ref{prop:ss3}) Let $X$ be the set $\{c\in\spp Q\st a\le ca\}$, and assume that $a\le \spp ba$. Then $\spp b\in X$, and thus assuming \ref{prop:ss4} we have
$\spp a=\bigwedge X\le\spp b$. The converse implication, namely that $\spp a\le \spp b$ implies $a\le\spp b a$, follows from (\ref{prop:s3a}).

(\ref{prop:s14}$\Rightarrow$\ref{prop:ss5}) Assume that \ref{prop:s14} holds and that $\spp a\le b 1$. Then $\spp a=\spp\spp a\le\spp(b 1)=\spp b$. Conversely, if $\spp a\le\spp b$ then
$\spp a\le\spp be\le\spp b 1=b1$.

(\ref{prop:ss5}$\Rightarrow$\ref{prop:s12}) We have $\spp(ab)\le abb^* a^*\le a1$, and thus assuming \ref{prop:ss5} we conclude $\spp(ab)\le\spp a$.

(\ref{prop:ss5}$\Rightarrow$\ref{prop:ss6}) If \ref{prop:ss5} holds then in the trivial equation $\spp a=\bigvee\{b\in\spp Q\st b\le \spp a\}$ we can substitute $b\le a 1$ for the equivalent condition $b\le \spp a$, thus obtaining $\spp a=\bigvee\{b\in\spp Q\st b\le a 1\}$.

(\ref{prop:ss6}$\Rightarrow$\ref{prop:ss5}) Let $X$ be the set $\{c\in\spp Q\st c\le b1\}$, and assume that
$\spp a\le b1$. Then $\spp a\in X$, and thus assuming \ref{prop:ss6} we have
$\spp a\le\V X=\spp b$. For the converse implication we note that $\spp a\le\spp b$ implies
$\spp a=\spp a e\le \spp b1=b1$.

(\ref{prop:s11}$\Rightarrow$\ref{prop:ss7}) The operation $(a,f)\mapsto a\cdot f=\spp(af)$ is a sup-lattice bimorphism, and $e$ acts trivially on $\spp Q$, since
$e\cdot f=\spp(e f)=\spp f=f$. Now let $a,b\in Q$ and $f\in\spp Q$. The associativity of the action follows from
\ref{prop:s11}: $(ab)\cdot f=\spp(a b f)=\spp(a\spp(b f))=a\cdot(b\cdot f)$. And the sup-lattice homomorphism $(-)1:\spp Q\to\rs(Q)$ preserves the action, for
$(a\cdot f)1=\spp(af)1=(af)1=a(f1)$.

(\ref{prop:ss7}$\Rightarrow$\ref{prop:s11}) Assume that $\spp Q$ is a module. Then \ref{prop:s11} follows from the associativity of the action: $\spp(ab)=\spp((ab)e)=(ab)\cdot e=a\cdot(b\cdot e)=\spp(a\spp b)$.
\qed
\end{proof}

\begin{definition}
A support is \emph{stable} if it satisfies the equivalent conditions of \ref{lemma:defstab}. A quantale equipped with a specified stable support is \emph{stably supported}.
\end{definition}

\begin{example}
All the examples of supports discussed so far are stable. A simple example of a supported quantale whose support is not stable is the four element unital involutive quantale
that, besides the elements $0$, $e$, and $1$, contains an element $a$ such that
\[\begin{array}{c}
a^2=a^*=a<e\;,\\
a 1=1\;.
\end{array}
\]
This quantale has a unique support, defined by $\spp a=a$, which is not stable: $\spp(a 1)=\spp 1=e\nleq a=\spp a$.
\end{example}

\begin{lemma}\label{prop:strongsupp}
Let $Q$ be a stably supported quantale.
\begin{enumerate}
\item Let $a,b\in Q$, and assume that the following three conditions hold:
\begin{eqnarray*}
b&\le& e\\
b&\le&a a^*\\
a&\le&b a\;.
\end{eqnarray*}
Then $b=\spp a$.
\item $\spp (a b)=a\spp b$ for all $a,b\in Q$ with $a\le e$.
\end{enumerate}
\end{lemma}

\begin{proof}
1. Assume $b\le e$. Then, from \ref{lemma:defstab}-\ref{prop:ss3}, the condition $a\le b a$ implies $\spp a\le b$. And the condition $b\le a a^*$ implies $b\le a 1$, which, by \ref{lemma:defstab}-\ref{prop:ss5}, is equivalent to
$b\le\spp a$.
Hence, if all the three conditions hold we conclude that $b=\spp a$.

2. If $a\le e$ we have $a\spp b\le e$, and thus
$a\spp b=\spp(a\spp b)=\spp(a b)$.
\qed
\end{proof}

\begin{theorem}\label{thm:fullsubcat}
\begin{enumerate}
\item If $Q$ has a stable support then that is the only support of $Q$, and the following equation holds:
\begin{equation}
\spp a = e\wedge a a^*\;.
\end{equation}
\item If $Q$ has a support and $K$ has a stable support then any homomorphism of unital involutive quantales from $Q$ to $K$ preserves the support. (In particular, the relational representations $Q\to\pwset{X\times X}$ of $Q$ are exactly the same as the support preserving relational representations.)
\end{enumerate}
\end{theorem}

\begin{proof}
1. Let $b=e\wedge a a^*$. Then by (\ref{def:s1}) and (\ref{def:s2}) we have $\spp a\le b$, and thus $a\le b a$, by \ref{lemma:defstab}-\ref{prop:ss3}. Hence, by \ref{prop:strongsupp}, we conclude that $b=\spp a$, which justifies the equation.

2. Let $h:Q\to K$ be a homomorphism of unital involutive quantales, let $a\in Q$, and let
$b=h(\spp a)$. Then we have:
\[
\begin{array}{l}
b = h(\spp a)\le h(e)=e\;,\\
b = h(\spp a)\le h(a a^*)=h(a) h(a)^*\;,\\
h(a)\le h(\spp a a)= h(\spp a) h(a)=b a\;.
\end{array}
\]
Hence, by \ref{prop:strongsupp} we conclude that $h(\spp a)=b=\spp h(a)$; that is, the support is preserved by $h$. \qed
\end{proof}

This theorem justifies the assertion that having a stable support is a property of a unital involutive quantale, rather than extra structure on it, and it motivates the following definition for the category of stably supported quantales (whose morphisms necessarily preserve the supports):

\begin{definition}
The \emph{category of stably supported quantales}, $\ssq$, is the full subcategory of the category of unital involutive quantales $\uiq$ whose objects are the stably supported quantales.
\end{definition}

In addition, the following theorem implies that each unital involutive quantale has an idempotent stably supported completion (idempotence meaning that any stably supported quantale is isomorphic to its completion, which is a consequence of $\ssq$ being a full subcategory):

\begin{theorem}\label{thm:fullreflective}
$\ssq$ is a full reflective subcategory of $\uiq$.
\end{theorem}

\begin{proof}
It is straightforward to see that the limits in $\ssq$ are calculated in $\uiq$, and thus the proof of the theorem will follow from verifying that the solution set condition of Freyd's adjoint functor theorem~\cite[Ch.\ V]{MacLane} holds.

In order to see this, first consider the category whose objects are the involutive monoids $M$ equipped with an additional operation $\spp:M\to M$, of which we require no special properties, and whose morphisms are the homomorphisms of involutive monoids that also preserve the operation $\spp$. Let us refer to such monoids as \emph{$\spp$-monoids}. From standard universal algebra it follows that there exist free
$\spp$-monoids.

Now let $X$ be a set and let us denote by $F(X)$ the corresponding free $\spp$-monoid. Let also $f:X\to K$ be a map, where $K$ is a supported quantale. Since the support makes $K$ an $\spp$-monoid, $f$ has a unique homomorphic extension $f':F(X)\to K$, of which there is then a unique join preserving extension $f'':\pwset{F(X)}\to K$. Besides being a sup-lattice, $\pwset{F(X)}$ is itself an $\spp$-monoid whose operations are computed pointwise from those of $F(X)$ (hence preserving joins in each variable). Furthermore, each of these operations is preserved by $f''$.

Now let $Q$ be a unital involutive quantale, and let
\[h:Q\to K\]
be a homomorphism of unital involutive quantales, where $K$ is stably supported. As above, there is a factorization:
\[
\xymatrix{Q\ar[rr]\ar[drr]_h&&\pwset{F(Q)}\ar[d]^{h''}\\
&&K.}
\]
Hence, the (necessarily stable) supported subquantale $S\subseteq K$ generated by the image $h(Q)$ is a surjective image of $\pwset{F(Q)}$, where the surjection $\pwset{F(X)}\to S$ is both a sup-lattice homomorphism and a homomorphism of $\spp$-monoids. This surjection determines an equivalence relation $\theta$ on $\pwset{F(Q)}$ such that:
\begin{itemize}
\item $\theta$ is a congruence with respect to joins, $\spp$, and the involutive monoid structure;
\item the injection of generators $Q \to\pwset{F(Q)}/\theta$ is a homomorphism of unital involutive quantales;
\item the quotient $\pwset{F(Q)}/\theta$ is stably supported.
\end{itemize}
Therefore we obtain a factorization
\[
\xymatrix{Q\ar[rr]^-{\eta_\theta}\ar[drr]_h&&\pwset{F(Q)}/\theta\ar[d]^{h'''}\\
&&K}
\]
in $\uiq$ (which implies that $h'''$ is in fact a homomorphism of supported quantales, by \ref{thm:fullsubcat}).
Since $K$ and $h$ have been chosen arbitrarily, the solution set condition now follows from the observation that the set of congruences which satisfy the above three conditions is small.

(One may also observe that the set of congruences is closed under intersections, and that the desired reflection is $\pwset{F(Q)}/\Theta$, where $\Theta$ is the least congruence.)
\qed
\end{proof}

\subsection{Partial units}\label{sec:partialunits}

\begin{definition}
Let $Q$ be a unital involutive quantale. A \emph{partial unit} in $Q$ is an element $a\in Q$ such that the following two conditions hold:
\begin{eqnarray*}
a a^*&\le& e\\
a^* a&\le& e\;.
\end{eqnarray*}
The set of partial units of $Q$ is denoted by $\ipi(Q)$.
\end{definition}

\begin{example}\label{exm:partialunits}
Let $X$ be a set, and $Q=\pwset{X\times X}$ the quantale of binary relations on $X$. Then $\ipi(Q)$ is the set $\ipi(X)$ of partial bijections on $X$. Hence, $\ipi(Q)$ is an inverse monoid and, as we shall see below, this is a consequence of the fact that $Q$ is a supported quantale.
More generally, if $Q=\pwset G$ for a discrete groupoid $G$, a partial unit is the same as a G-set in the sense of \cite{Renault} (see section \ref{sec:prelim.invsems}), and thus $\ipi(Q)=\ipi(G)$.
\end{example}

\begin{lemma}
Let $Q$ be a unital involutive quantale. Then $\ipi(Q)$ is an involutive submonoid of $Q$.
\end{lemma}

\begin{proof}
The set $\ipi(Q)$ is clearly closed under involution, and $e\in\ipi(Q)$. It is also closed under multiplication, for if $a$ and $b$ are partial units then $(a b)(a b)^*=a b b^* a^*\le a e a^*=a a^*\le e$, and in the same way $(a b)^*(a b)\le e$. Hence, $\ipi(Q)$ is an involutive submonoid of $Q$. \qed
\end{proof}

\begin{lemma}\label{prop:ipi}
Let $Q$ be a supported quantale, and let $a\in \ipi(Q)$. Then,
\begin{enumerate}
\item $\spp a=a a^*$,
\item $a=a a^* a$,
\item $a^2=a$ if and only if $a\le e$,
\item $b\le a$ if and only if $b=\spp b a$, for all $b\in Q$.
\end{enumerate}
\end{lemma}

\begin{proof}
1. We have $\spp a\le a a^*\le\spp a a a^*\le\spp a e=\spp a$.

2. This is a consequence of the above and of the equality $a=\spp a a$.

3. If $a\le e$ then $a$ is an idempotent because $\spp Q$ is a frame. Now assume that $a^2=a$. Then
$a a^*= a a a^*\le a e=a$.
Hence, $a a^*\le a$, and, since $a a^*$ is self-adjoint, also $a a^*\le a^*$. Finally, from here we conclude that $a\le e$ because $a=a a^* a\le a^* a\le e$.

4. Let $b\in Q$ such that $b\le a$ (in fact, then $b\in\ipi(Q)$). Then,
\[\spp b a\le b b^* a\le b a^* a\le b e=b\le\spp b b\le\spp b a\;.\]
This shows that $b\le a$ implies $b=\spp b a$. The converse is trivial. \qed
\end{proof}

In the following theorem the category of inverse quantales should be naturally assumed to be the category whose objects are the supported quantales and whose morphisms are the unital and involutive homomorphisms that preserve the supports, although the theorem is true even if we consider as morphisms all the homomorphisms of unital involutive quantales (in any case the distinction disappears once we restrict to stable supports).

\begin{theorem}\label{thm:ipiQ}
\begin{enumerate}
\item
Let $Q$ be a supported quantale. Then $\ipi(Q)$ is an inverse monoid whose natural order coincides with the order inherited from $Q$, and whose set of idempotents $E(\ipi(Q))$ coincides with $\spp Q$.
\item The assignment $Q\mapsto\ipi(Q)$ extends to a functor $\ipi$ from the category of supported quantales to the category of inverse monoids $\imcat$.
\end{enumerate}
\end{theorem}

\begin{proof}
1. $\ipi(Q)$ is an involutive submonoid of $Q$, and in particular it is a regular monoid because for each partial unit $a$ we have both $a a^* a=a$ and $a^* a a^*=a^*$. Hence, in order to have an inverse monoid it suffices to show that all the idempotents commute, and this follows from \ref{prop:ipi}-3, which implies that the set of idempotents of $\ipi(Q)$ is the same as $\spp Q$, which is a frame. Furthermore, the natural order of $\ipi(Q)$ is defined by $a\le b\Leftrightarrow\spp a b= a$, and thus it coincides with the order of $Q$, by \ref{prop:ipi}-4.

2. $\ipi$ is a functor because if $h:Q\to K$ is any homomorphism of unital involutive quantales and $a\in Q$ is a partial unit then $h(a)$ is a partial unit: $h(a) h(a)^*=h(a a^*)\le h(e)=e$, and, similarly, $h(a)^* h(a)\le e$. Hence, $h$ restricts to a homomorphism of monoids $\ipi(Q)\to\ipi(K)$.
\qed
\end{proof}

\subsection{Inverse quantales}\label{sec:inversequantales}

We shall be particularly interested in supported quantales with the additional property that each element is a join of partial units. As we shall see, such quantales are necessarily stably supported.

\begin{definition}
By an \emph{inverse quantale} will be meant a supported quantale $Q$ such that every element $a\in Q$ is a join of partial units:
\[a=\V\{b\in\ipi(Q)\st b\le a\}\;.\]
The \emph{category of inverse quantales}, $\iqcat$, is the full subcategory of $\uiq$ whose objects are the inverse quantales.
\end{definition}

From \ref{exm:partialunits} it follows that $\pwset G$ is an inverse quantale for any discrete groupoid $G$. In particular, the quantale of binary relations $\pwset{X\times X}$ on a set $X$ is an inverse quantale.

The following result provides an alternative definition of the concept of inverse quantale.

\begin{lemma}\label{thm:inversequantale}
Let $Q$ be a unital involutive quantale such that every element in $Q$ is a join of partial units. Then $Q$ is a supported quantale (hence, an inverse quantale) if and only if the following conditions hold:
\begin{enumerate}
\item $a\le aa^*a$ for all $a\in Q$ (equivalently, $a=aa^*a$ for all $a\in\ipi(Q)$);
\item the operation $\spp:Q\to Q$ defined by
\[\spp a=\V\{bb^*\st b\in\ipi(Q)\textrm{ and }b\le a\}\]
is a sup-lattice homomorphism.
\end{enumerate}
When these conditions hold the operation $\spp$ is the support.
\end{lemma}

\begin{proof}
If $Q$ is a supported quantale (equivalently, an inverse quantale) then the condition $a\le aa^*a$ holds, and the above formula for the support follows from the fact that for each partial unit $b\in\ipi(Q)$ we have $\spp b=bb^*$ (\cf\ \ref{prop:ipi}):
\begin{eqnarray*}
\spp a&=&\spp\left({\V\{b\st b\in\ipi(Q),\ b\le a\}}\right)=\V\{\spp b\st b\in\ipi(Q),\ b\le a\}\\
&=&\V\{bb^*\st b\in\ipi(Q),\ b\le a\}\;.
\end{eqnarray*}
Conversely, let $Q$ be a unital involutive quantale each of whose elements is a join of partial units, in addition satisfying $a\le aa^*a$ for all $a\in Q$, and such that the operation $\spp$ defined by the above formula preserves joins. This operation satisfies the conditions (\ref{def:s1}), (\ref{def:s2}), and (\ref{def:s3}), and thus it is a support:
\begin{itemize}
\item
by definition of partial unit we have $bb^*\le e$ for all $b\in\ipi(Q)$, and thus $\spp a\le e$;
\item we have
$\spp a\le aa^*$ because $bb^*\le aa^*$ for all $b\le a$;
\item we have
\begin{eqnarray*}
\spp a a&=&\V\{bb^* a\st b\in\ipi(Q)\textrm{ and }b\le a\}\\
&\ge&\V\{bb^*b\st b\in\ipi(Q)\textrm{ and }b\le a\}\\
&=&\V\{b\st b\in\ipi(Q)\textrm{ and }b\le a\}=a\;.\qed
\end{eqnarray*}
\end{itemize}
\end{proof}

\begin{example}
The following commutative unital quantale $Q$, with trivial involution, shows that the second condition of the above lemma above is not redundant:
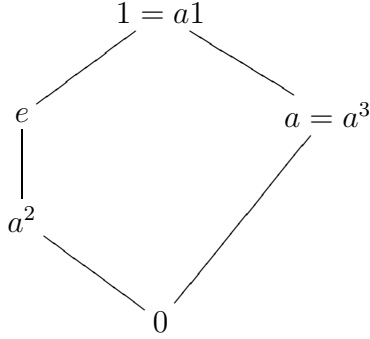
\[
\xymatrix{
&1=a1\ar@{-}[dr]\ar@{-}[dl]\\
e\ar@{-}[d]&&a=a^3\ar@{-}[ddl]\\
a^2\ar@{-}[dr]\\
&0
}
\]
In this example, which clearly satisfies the condition $b\le bb^*b$ for all $b$, every element is a join of partial units, but the operation $\spp$ defined by
\[\spp a=\V\{bb^*\st b\in\ipi(Q)\textrm{ and }b\le a\}\]
does not preserve joins, for $\spp a=aa^*=a^2$ and thus $\spp(a^2)\vee\spp a=\spp\spp a\vee\spp a=\spp a=a^2$, whereas $\spp(a^2\vee a)=\spp 1=e\neq a^2$.
\end{example}

\begin{lemma}\label{thm:invimpliesss}
Any inverse quantale is stably supported (hence, $\iqcat$ is a full subcategory of $\ssq$).
\end{lemma}

\begin{proof}
Let $Q$ be an inverse quantale, and let $a,b\in Q$, such that
\begin{eqnarray*}
a&=&\V_i s_i\\
b&=&\V_j t_j\;,
\end{eqnarray*}
where the $s_i$ and $t_j$ are partial units.
Then,
\begin{eqnarray*}
\spp (ab)&=&\spp\left({\V_{ij} s_i t_j}\right)
=\V_{ij}\spp(s_i t_j)
=\V_{ij} s_i t_j(s_i t_j)^*
=\V_{ij} s_i t_j t_j^* s_i^*\\
&\le&\V_i s_i s_i^*
=\V_i \spp(s_i)
=\spp\left({\V_i s_i}\right)
=\spp a\;. \qed
\end{eqnarray*}
\end{proof}

The converse is not true; that is, not every stably supported quantale is an inverse quantale (see \ref{fuzzyexample} in section \ref{sec:examples}).

\begin{lemma}\label{lemma:invqfrominvs}
Let $S$ be an inverse semigroup. The set $\lc(S)$ of subsets of $S$ that are downwards closed in the natural order of $S$ is an inverse quantale. The unit is the set of idempotents $E(S)$ (which, if $S$ is a monoid, is just $\downsegment e$), multiplication is computed pointwise,
\[X Y=\{xy\st x\in X,\ y\in Y\}\;,\]
the involution is pointwise inversion, $X^*=X^{-1}$, and the (necessarily unique) support is given by the formula
\[\spp X=\{xx^{-1}\st x\in X\}\;.\]
\end{lemma}

\begin{proof}
Consider an inverse semigroup $S$. It is straightforward to verify that the sup-lattice $\lc(S)$ is an involutive quantale, with the multiplication and involution defined above. In particular, for multiplication this means that downwards closed sets are closed under pointwise multiplication (this would not be true for an arbitrary partially ordered involutive semigroup, for which downwards closure would be required after taking the pointwise multiplication), which is a consequence of the fact that we are dealing with the natural order of an inverse semigroup: if $z\le w\in XY$ for $X,Y\in\lc(S)$, then $w=xy$ with $x\in X$ and $y\in Y$, and $z=zz^{-1}w=zz^{-1}(xy)=(zz^{-1}x)y\in XY$ because $zz^{-1}x\le x\in X$, and $X$ is downwards closed. For the involution it is similar, but more immediate.

In order to see that the set of idempotents $E(S)$ is the multiplicative unit, consider $X\in\lc(S)$. The set $E(S)$ contains the idempotent $x^{-1}x$ for each $x\in X$, and thus the pointwise product
$XE(S)$ contains all the elements of the form $xx^{-1}x=x$. Hence, $X\subseteq XE(S)$. The other elements of $XE(S)$ are of the form $xy^{-1}y$, with $x\in X$, and we have $xy^{-1}y\le x$ in the natural order of $S$, implying that $xy^{-1}y\in X$ because $X$ is downwards closed. Hence, $X E(S)=X$. Similarly we show that $E(S) X=X$.

In order to see that $\lc(S)$ is an inverse quantale we apply \ref{thm:inversequantale}. First, each $X\in\lc(S)$ is of course a union of partial units:
\[X=\bigcup\{\downsegment x\st x\in X\}\;.\]
Secondly, let $X\in\lc(S)$. Then
\[XX^*X=\{xy^{-1}z\st x,y,z\in X\}\supseteq\{xx^{-1}x\st x\in X\}=\{x\st x\in X\}=X\;.\]
To conclude, we show that the operation $\spp:\lc(S)\to\lc(S)$ defined by
\[\spp X=\bigcup\{UU^*\st U\in\ipi(\lc(S))\textrm{ and }U\subseteq X\}\]
coincides with the (clearly join-preserving) operation
\[\spp X=\{xx^{-1}\st x\in X\}\]
of the statement of this lemma, thus showing that $\lc(S)$ is an inverse quantale.
If $U\in\ipi(\lc(S))$ then $UU^*\subseteq E(S)$, \ie, $xy^{-1}\in E(S)$ for all $x,y\in U$. But then
\[xy^{-1}=xy^{-1}(xy^{-1})^{-1}=xy^{-1}yx^{-1}=xy^{-1}yy^{-1}yx^{-1}=(xy^{-1}y)(xy^{-1}y)^{-1}\;,\]
so we see that $xy^{-1}$ is of the form $zz^{-1}$ with $z=xy^{-1}y\le x$, where $z\in U$ because $U$ is downwards closed. We therefore conclude that $UU^*$ coincides with the set $\{xx^{-1}\st x\in U\}$, and from this the required formula for $\spp X$ follows. \qed
\end{proof}

\begin{theorem}\label{thm:LS}
The functor $\ipi$, restricted to the category of stably supported quantales $\ssq$, has a left adjoint from the category of inverse monoids $\imcat$ to $\ssq$, which to each inverse monoid $S$ assigns the quantale $\lc(S)$.
\end{theorem}

\begin{proof}
First we remark that the embedding $S\to\lc(S)$ actually defines a homomorphism of monoids
$S\to\ipi(\lc(S))$, which provides the unit of the adjunction.
Now let $Q$ be a stably supported quantale, and
$h:S\mapsto\ipi(Q)$ a homomorphism of monoids. Then $h$ preserves the natural order, and thus it defines a homomorphism of ordered involutive monoids $S\to Q$ because, as we have seen, the natural order of $\ipi(Q)$ is just the order of $Q$ restricted to $\ipi(Q)$. It follows that $h$ extends (uniquely) to a homomorphism of unital involutive quantales $\overline h:\lc(S)\to Q$, namely the sup-lattice extension
$\overline h(U)=\V h(U)$ of \ref{prop:supextension}. \qed
\end{proof}

\subsection{Enveloping quantales}\label{sec:envelopes}

We have seen how to obtain an inverse quantale from an arbitrary inverse semigroup. In the case of an abstract complete pseudogroup there is another useful inverse quantale $\lcc(S)$, here referred to as the \emph{enveloping quantale} of $S$, which takes into account the joins that exist in $S$. From here until the end of this section $S$ will be a fixed but arbitrary abstract complete pseudogroup.

\begin{definition}
By a \emph{compatibly closed ideal} of $S$ is meant a downwards closed set (possibly empty) which is closed under the formation of joins of compatible sets. The set of compatibly closed ideals of $S$ is denoted by $\lcc(S)$.
\end{definition}

\begin{lemma}
$\lcc(S)$ is a quotient of $\lc(S)$ both as a frame and as a unital involutive quantale, and it is an inverse quantale.
\end{lemma}

\begin{proof}
Let $j:\lc(S)\to\lc(S)$ be the closure operator that to each downwards closed set $U\subseteq S$ assigns the least compatibly closed ideal that contains $U$. First we remark that $j$ is explicitly defined by
\[j(U)=\left\{{\V X\st X\subseteq U,\ X\textrm{ is compatible}}\right\}\;.\]
In order to see this, let $x\le y\in j(U)$. Then $y$ is of the form $\V Y$ for some compatible set $Y\subseteq U$, and thus $x=x\wedge\V Y=\V (x\wedge Y)$, where the set $x\wedge Y=\{x\wedge y\st y\in Y\}$ is of course compatible. Hence, we have $x\in j(U)$, showing that $j(U)$ is downwards closed. Now let $Z\subseteq j(U)$ be a compatible set. Each element $z\in Z$ is of the form $\V U_z$ for some compatible set $U_z\subseteq U$, and the fact that $Z$ is compatible implies that the set $Z'=\bigcup_{z\in Z} U_z$ is compatible. But we also have $\V Z=\V Z'$, and thus $z\in j(U)$, showing that $j(U)$ is closed under the formation of joins of compatible sets. It is thus a compatible ideal, clearly the smallest one containing $U$.

Now we shall show that $j$ is both a frame nucleus and a nucleus of involutive quantales, thus proving that $\lcc(S)$, which coincides with the quotient of $\lc(S)$ obtained as the set of fixed-points of $S$, is a unital involutive quantale, a frame, and a quotient of $\lc(S)$.

Let $I,J\in\lc(S)$. Let $X\subseteq I$ and $Y\subseteq J$ be compatible sets such that $\V X=\V Y$. Let us denote this join by $z$. We have $z\in j(I)\cap j(J)$, and all the elements of $j(I)\cap j(J)$ can be obtained in the same way. Now define the set
\[Z=\{x\wedge y\st x\in X,\ y\in Y\}\;.\]
We have $Z\subseteq I\cap J$, and
\[\V Z=\V_{x\in X}\V_{y\in Y}(x\wedge y)=\V_{x\in X}\left({x\wedge\V Y }\right)=\V_{x\in X}(x\wedge z)=\left({\V X}\right)\wedge z=z\;.\]
Hence, $z\in j(I\cap J)$, and we conclude that $j(I)\cap j(J)\subseteq j(I\cap J)$, \ie, $j$ is a nucleus of frames.

Let again $I,J\in\lc(S)$. Consider an arbitrary element of $j(I)j(J)$, which is necessarily of the form $xy$ with $x=\V X$ and $y=\V Y$, where $X\subseteq I$ and $Y\subseteq J$ are compatible sets. We shall show that $xy\in j(XY)$, hence proving that $j$ is a nucleus of quantales. First, we remark that $XY$ is a compatible set, since it is bounded above by $xy$. Its join $\V(XY)$ coincides with $xy$, due to infinite distributivity, and thus $xy\in j(XY)$. Clearly, $j$ preserves the involution of the quantales, and thus $j$ is a nucleus of involutive quantales. The involution of $\lcc(S)$ is, similarly to that of $\lc(S)$, given by pointwise inversion.

Finally, $j$ is also a nucleus with respect to the support of $\lc(S)$ because the support is, since $\lc(S)$ is stably supported, expressed in terms of frame and unital involutive quantale operations:
\[\spp U = UU^*\cap E(S)\;.\]
Hence, the conclusion that $\lcc(S)$ is stably supported follows, and it is obviously an inverse quantale because it has a basis consisting of the principal ideals $\downsegment s$ with $s\in S$, which are partial units of $\lcc(S)$. \qed
\end{proof}

\begin{example}
From the results about groupoids and quantales later in this paper it will follow (but it can also be verified directly) that if $G$ is a discrete groupoid then we have an isomorphism
$\lcc(\ipi(G))\cong\pwset G$. In particular, the enveloping quantale $\lcc(\ipi(X))$ of the symmetric inverse monoid $\ipi(X)$ of a set $X$ is isomorphic to the quantale of binary relations $\pwset{X\times X}$ on $X$.
\end{example}

\begin{theorem}
$\lcc(S)$ is the quotient of $\lc(S)$ (in the category of stably supported quantales $\ssq$) determined by the condition that joins of $E(S)$ should be preserved by the injection of generators $S\to\lcc(S)$.
\end{theorem}

\begin{proof}
This is a corollary of \ref{prop:abspg}, from which it follows that the homomorphisms of abstract complete pseudogroups (\ie, those monoid homomorphisms that preserve joins of compatible subsets)
\[h:S\to\ipi(Q)\;,\]
where $Q$ is any stably supported quantale, are exactly the monoid homomorphisms that preserve just the joins of sets of idempotents. Hence, the universal properties corresponding to preservation of joins of idempotents, on one hand, and to preservation of arbitrary joins, on the other, are the same. \qed
\end{proof}

The universal properties possessed by enveloping quantales are by now essentially obvious, once one takes into account the analogous properties for $\lc(S)$. We shall provide an explicit description of them.

Let $S$ be an abstract complete pseudogroup, and let $Q$ be a stably supported quantale. The monoid of partial units $\ipi(Q)$ is an abstract complete pseudogroup, and if $h:S\to\ipi(Q)$ is a homomorphism of abstract complete pseudogroups there is a unique homomorphism of unital involutive quantales $\bar h:\lcc(S)\to Q$ such that the following diagram commutes,
\[\xymatrix{S\ar[rr]^-{s\mapsto\downsegment s}\ar[drr]_h&&\lcc(S)\ar[d]^{\bar h}\\
&& Q}\]
where $\bar h$ is explicitly defined by
\[\bar h(U)=\V\{h(s)\in Q\st s\in U\}\;.\]
In other words, we have:

\begin{corollary}\label{cor:adjunction}
$\lcc$ defines a functor from $\apgcat$ to $\ssq$, which is left adjoint to the functor $\ipi:\ssq\to\apgcat$.
\end{corollary}

A consequence of this is also, for $S$ an inverse semigroup, that $\lc(S)\cong\lcc(C(S))$, where $C(S)$ is the completion of $S$ in the sense of \cite[Section 1.4]{Lawson} (but including the join of the empty set), and in fact we have $C(S)\cong\ipi(\lc(S))$.

\section{Quantal frames}

We have already remarked that the topologies of certain topological groupoids are quantales. Besides this they are also frames, of course, hence suggesting the following definition:

\begin{definition}\label{def:quantalframe}
By a \emph{quantal frame} is meant a quantale $Q$ such that for all $a,b_i\in Q$ the following distributivity property holds:
\[a\wedge\V_i b_i=\V_i a\wedge b_i\;.\]
\end{definition}

In this section we shall see how (localic) groupoids or at least categories can be obtained from suitable quantal frames. The main result (theorem \ref{inversiontheorem}) states simply that the localic \'{e}tale groupoids correspond bijectively, up to isomorphisms, to the inverse quantales that are also quantal frames. To a large extent this is a consequence of a bijective correspondence between these quantales and the abstract complete pseudogroups that furthermore is part of an equivalence of categories (theorem \ref{thm:equivapgiqf}). At the end, in section \ref{sec:examples}, we provide examples whose purpose is to separate all the classes of quantales considered so far.

\subsection{Stable quantal frames}\label{sec:sqf}

We begin by considering quantal frames of which nothing is required except that their underlying quantales should be stably supported. This condition can be expressed equivalently as follows:

\begin{definition}\label{def:stablequantale}
By a \emph{stable quantal frame} will be meant a unital involutive quantal frame satisfying the following additional conditions:
\begin{eqnarray*}
a1\wedge e&\le& aa^*\\
a&\le&(a1\wedge e)a\;.
\end{eqnarray*}
[Equivalently, satisfying the equations $a1\wedge e= aa^*\wedge e$ and
$a=(a1\wedge e)a$.]
\end{definition}

Let $Q$ be a stable quantal frame.
The sup-lattice inclusion
\[\upsilon:\spp Q\to Q\;.\]
has a right adjoint given by
\[\upsilon_*(a)=a\wedge e\;,\]
which preserves arbitrary joins and is the inverse image frame homomorphism of an (obviously) open locale map
\[u:G_0\to G_1\]
whose direct image is $u_!=\upsilon$,
where $G_0$ and $G_1$ are defined by the conditions
\[\opens(G_0)=\spp Q\textrm{ and }\opens(G_1)=Q\;.\]

Now consider the sup-lattice homomorphism
$\delta:Q\to\spp Q$
defined by
\[
\delta(a)=a1\wedge e\;.
\]
(This is just the support $\spp$ with codomain restricted to $\spp Q$.)

\begin{lemma}
$\delta$ is the direct image $d_!$ of an open map $d:G_1\to G_0$.
\end{lemma}

\begin{proof}
Consider the map $\spp Q\to Q$ given by $a\mapsto a1$. By the properties of stably supported quantales this is an isomorphism $\spp Q\to \rs(Q)$ followed by the inclusion $\rs(Q)\to Q$. Hence, it is a frame homomorphism and it defines a map of locales $d:G_1\to G_0$, which furthermore is semiopen with $\delta=d_!$ because $d^*$ is the right adjoint $\delta_*$ of $\delta$:
\[
\begin{array}{c}
\delta(d^*(a))=\spp(a1)=\spp a=a\textrm{ for all }a\le e\\
d^*(\delta(a))=\spp a 1=a1\ge a\textrm{ for all }a\in Q\;.
\end{array}
\]
In order to see that $d$ is open we check the Frobenius reciprocity condition. Let $a,b\in Q$, with $b\le e$. Then
\[
\begin{array}{rcll}
d_!(d^*(b)\wedge a)&=&\spp(b1\wedge a)\\
&=&\spp(b a)&\textrm{(By \ref{lemma:defstab}-\ref{prop:s15a} and $\spp b=b$.)}\\
&=&b\spp a&\textrm{(By \ref{prop:strongsupp}-2.)}\\
&=&b\wedge \spp a&\textrm{($\spp Q$ is a locale.)}\\
&=&b\wedge d_!(a)\;. \qed
\end{array}
\]
\end{proof}

The involution of $Q$ is a frame isomorphism that defines a locale map $i:G_1\to G_1$ by the condition $i^*(a)=a^*$, and thus we have $i\circ i=\ident$ and $i_!= i^*$.
Our aim is that ultimately $d$ should be the domain map of a groupoid, and we obtain a candidate for an open range map
\[r:G_1\to G_0\]
just by defining $r=d\circ i$.
These maps satisfy the appropriate relations:

\begin{lemma}\label{lemma:graphconditions}
Consider the locale maps
\[
\xymatrix{G_1\ar@(ld,lu)^i\ar@<1.5ex>[rr]^d\ar@<-1.5ex>[rr]_r&&G_0\ar[ll]|u}
\]
as defined above. We have
\begin{eqnarray*}
d\circ u&=&\ident\\
r\circ u&=&\ident\\
d\circ i&=&r\\
r\circ i&=&d\;.\\
\end{eqnarray*}
\end{lemma}

\begin{proof}
The first condition is equivalent to $u^* \circ d^*=\ident$, which holds because for all $a\le e$ we have $\spp a=a$, and thus
\[u^*(d^*(a))=a1\wedge e=\spp a=a\;.\]
Similarly, the second condition is true: for all $a\le e$ we have $a=a^*$, and thus
\[u^*(r^*(a))=\spp(a^*)=\spp a=a\;.\]
The third condition is the definition of $r$, and the fourth follows from this because $i\circ i=\ident$. \qed
\end{proof}

So far we have obtained from the stable quantal frame $Q$ a localic graph $G$ that is equipped with an involution $i$, and whose maps $d$, $r$, and $u$, are all open. There is additional structure on $G$, consisting of a certain kind of multiplication defined on the ``locale of composable pairs of edges'' $G_1\times_{G_0}G_1$, although not (yet) necessarily the multiplication of a groupoid or even of a category. In order to see this let us first notice that the frame $\opens(G_0)=\spp Q$ is, as a quantale, a unital subquantale of $Q$ because meet in $\spp Q$ coincides with multiplication in $Q$ (see \ref{prop:supps}). Hence, there are two immediate ways in which $Q$ is a module over $\spp Q$: multiplication on the left defines an action of $\spp Q$ on $Q$, and multiplication on the right defines another. We shall regard $Q$ as an $\spp Q$-$\spp Q$-bimodule with respect to these two actions, namely letting the left (resp.\ right) action be left (resp.\ right) multiplication (in fact each action makes $Q$ both a right and a left $\spp Q$-module because $\spp Q$ is a commutative quantale, but this is irrelevant). We shall denote the corresponding tensor product over $\spp Q$ by $Q\otimes_{\spp Q}Q$. This coincides (fortunately for our notation) with $\opens(G_1\times_{G_0}G_1)$:

\begin{lemma}
The tensor product $Q\otimes_{\spp Q}Q$ coincides with the pushout of the homomorphisms $d^*$ and $r^*$.
\end{lemma}

\begin{proof}
This is equivalent to showing, for all $a,b,c\in Q$, with $a\le e$, that the equality
\[(b\wedge r^*(a))\otimes c=b\otimes (d^*(a)\wedge c)\]
(see the end of section \ref{sec:toplocgrpds}) is equivalent to
\[ba\otimes c=b\otimes ac\;,\]
which is immediate from \ref{lemma:defstab}-\ref{prop:s15a}:
\[d^*(a)\wedge c=a1\wedge c=\spp a c=ac\;.\]
[For $r^*$ it is analogous, using the obvious dual of \ref{lemma:defstab}-\ref{prop:s15a}.] \qed
\end{proof}

Now we notice that the quantale multiplication
$Q\otimes Q\to Q$
factors through the above pushout because it is associative and thus in particular it respects the relations $ba\otimes c=b\otimes ac$ for $a\le e$ that determine the quotient
$Q\otimes Q\to Q\otimes_{\spp Q}Q$. This enables us to make the following definition:

\begin{definition}\label{def:quantalmultiplication}
Let $Q$ be a stable quantal frame, and let
\[G=\xymatrix{G_1\ar@(ld,lu)^i\ar@<1.5ex>[rr]^d\ar@<-1.5ex>[rr]_r&&G_0\ar[ll]|u}\]
be the corresponding involutive localic graph. The \emph{quantal multiplication induced by $Q$ on $G$} is the sup-lattice homomorphism
\[\mu:\opens(G_1\ptimes{G_0} G_1)\to \opens(G_1)\]
which is defined by, for all $a,b\in Q$,
\[\mu(a\otimes b)=ab\;.\]
\end{definition}

\subsection{Multiplicative quantal frames}

Let us now examine a condition under which the localic graph associated to a stable quantal frame has the additional structure of a localic category.

\begin{definition}
By a \emph{multiplicative quantal frame} will be meant a stable quantal frame for which the right adjoint $\mu_*$ of the quantal multiplication $\mu:Q\otimes_{\spp Q}Q\to Q$ preserves arbitrary joins.
\end{definition}

The multiplicativity condition is that under which $\mu$ is the direct image $m_!$ of a (semiopen) locale map $m$, which then gives us a category, as the following theorem shows.

\begin{theorem}\label{theorem:localiccategory}
Let $Q$ be a multiplicative quantal frame. Then the locale map
\[m:G_1\ptimes{G_0}G_1\to G_1\]
which is defined by $m^*=\mu_*$ (equivalently, $m_!=\mu$),
together with the maps $d$, $r$, and $u$, defines a localic category.
\end{theorem}

\begin{proof}
In \ref{lemma:graphconditions} we have obtained many of the needed conditions. The only ones missing are the unit laws, the associativity of $m$, and those that specify the domain and range of a product of two arrows:
\begin{eqnarray}
d\circ m &=& d\circ\pi_1 \label{simp1}\\
r\circ m&=& r\circ\pi_2\;. \label{simp2}
\end{eqnarray}
We shall begin by proving these.
In fact we shall prove (\ref{simp1}) only, as (\ref{simp2}) is analogous. We have to show, for frame homomorphisms, that $m^*\circ d^*=\pi_1^*\circ d^*$. In order to do this we shall prove that $d_!\circ m_!$ is left adjoint to $\pi_1^*\circ d^*$ (this identifies $m^*\circ d^*$ and $\pi_1^*\circ d^*$ because adjoints between partial orders are uniquely determined), in order to take advantage of the following simple formulas:
\begin{eqnarray*}
d_!(m_!(a\otimes b))&=&\spp(ab)\\
\pi_1^*(d^*(a))&=&a1\otimes 1\hspace*{1cm}(a\le e)
\end{eqnarray*}
We prove that the adjunction exists by proving the following two inequalities (resp.\ the co-unit and the unit of the adjunction):
\begin{eqnarray*}
d_!(m_!(\pi_1^*(d^*(a))))&\le& a\textrm{ for all }a\le e\\
\pi_1^*(d^*(d_!(m_!(a\otimes b))))&\ge& a\otimes b\textrm{ for all }a,b\in Q\;.
\end{eqnarray*}
Let us prove the first inequality. Consider $a\le e$. Then
\[d_!(m_!(\pi_1^*(d^*(a))))=d_!(m_!(a1\otimes 1))=\spp(a11)=\spp(a1)=\spp a=a\;.\]
Now let us prove the second inequality. Let $a,b\in Q$. Then
\begin{eqnarray*}
\pi_1^*(d^*(d_!(m_!(a\otimes b))))&=&\pi_1^*(d^*(\spp(ab)))=\spp(ab)1\otimes 1=\\
&=&ab1\otimes 1=a\spp b 1\otimes 1\ge a\spp b\otimes 1=\\
&=&a\otimes\spp b 1=a\otimes b1\ge a\otimes b\;.
\end{eqnarray*}
Hence, (\ref{simp1}) holds, and for analogous reasons so does (\ref{simp2}).

Let us now prove the unit laws, which state that the following diagram is commutative:
\begin{equation}
\label{unitdiagram}
\vcenter{\xymatrix{
G_0\ptimes{G_0}G_1\ar[rr]^{u\times\ident}&&G_1\ptimes{G_0}G_1\ar[d]_m&&
G_1\ptimes{G_0}G_0\ar[ll]_{\ident\times u}\\
G_1\ar[u]^-{\langle d,\ident\rangle}\ar@{=}[rr]&&G_1&&G_1\ar[u]_-{\langle\ident,r\rangle}\ar@{=}[ll]
}}
\end{equation}
First we remark that, since the frame homomorphism $u^*$ has a left adjoint $u_!$, the maps $u\times \ident$ and $\ident\times u$ are semiopen, with
\begin{eqnarray*}
(u\times\ident)_!&=&u_!\otimes\ident\\
(\ident\times u)_!&=&\ident\otimes u_!\;,
\end{eqnarray*}
because the operations $\ident\otimes -$ and $-\otimes\ident$ are functorial and thus preserve the conditions $u_!\circ u^*\le\ident$ and $u^*\circ u_!\ge\ident$ that define the adjunction $u_!\dashv u^*$.
Secondly, the maps $\langle d,\ident\rangle$ and $\langle\ident,r\rangle$ are isomorphisms whose inverses are, respectively, the projections $\pi_2:G_0\times_{G_0}G_1\to G_1$ and $\pi_1:G_1\times_{G_0}G_0\to G_1$. Hence, in particular, these maps are semiopen, and we have
\begin{eqnarray*}
\langle d,\ident\rangle_!&=&\pi_2^*\\
\langle\ident,r\rangle_!&=&\pi_1^*\;.
\end{eqnarray*}
Hence, since $m^*$, too, has a left adjoint $m_!$, we conclude that the commutativity of
(\ref{unitdiagram}) is equivalent to that of the following diagram:
\begin{equation}
\label{unitdiagram2}
\vcenter{\xymatrix {
\opens(G_0\ptimes{G_0}G_1)\ar[rr]^{u_!\otimes\ident}&&\opens(G_1\ptimes{G_0}G_1)\ar[d]_{m_!}&&
\opens(G_1\ptimes{G_0}G_0)\ar[ll]_{\ident\otimes u_!}\\
\opens(G_1)\ar[u]^-{\pi_2^*}\ar@{=}[rr]&&\opens(G_1)&&\opens(G_1)\ar[u]_-{\pi_1^*}\ar@{=}[ll]
}}
\end{equation}
Now the commutativity of the left square of (\ref{unitdiagram2}), \ie, the condition
\[m_!\circ (u_!\otimes\ident)\circ\pi_2^*=\ident\;,\]
follows from the fact that $Q$ is a unital quantale with $e=u_!(1_{G_0})$, since for all $a\in Q$ we obtain
\[(m_!\circ (u_!\otimes\ident)\circ\pi_2^*)(a)=m_!\circ(u_!\otimes\ident)(1_{G_0}\otimes a)=m_!(u_!(1_{G_0})\otimes a)=ea=a\;.\]
Similarly, the right square of (\ref{unitdiagram2}) follows from $ae=a$.

Finally, we shall prove that the multiplication map $m$ is associative.
Let $m':\opens(G_1)\otimes\opens(G_1)\to\opens(G_1)$ be the quantale multiplication,
\[m'=m_!\circ q\;,\]
where $q$ is is the quotient homomorphism
\[q:\opens(G_1)\otimes\opens(G_1)\to\opens(G_1)\ootimes{\opens(G_0)}\opens(G_1)\;.\]
It is clear that the associativity of $m'$, which is equivalent to the commutativity of
\[
\xymatrix{
(\opens(G_1)\otimes\opens(G_1))\otimes\opens(G_1)\ar[rrrr]^-{m'\otimes \ident}\ar[d]_{\cong}&&&&\opens(G_1)\otimes\opens(G_1)\ar[dd]^{m'}\\
\opens(G_1)\otimes(\opens(G_1)\otimes\opens(G_1))\ar[d]_{\ident\otimes m'}\\
\opens(G_1)\otimes\opens(G_1)\ar[rrrr]^-{m'}&&&&\opens(G_1)\;,
}
\]
implies (in fact it is equivalent to) the commutativity of
\[
\xymatrix{
(\opens(G_1)\ootimes{\opens(G_0)}\opens(G_1))\ootimes{\opens(G_0)}\opens(G_1)\ar[rrrr]^-{m_!\otimes \ident}\ar[d]_{\cong}&&&&\opens(G_1)\ootimes{\opens(G_0)}\opens(G_1)\ar[dd]^{m_!}\\
\opens(G_1)\ootimes{\opens(G_0)}(\opens(G_1)\ootimes{\opens(G_0)}\opens(G_1))\ar[d]_{\ident\otimes m_!}\\
\opens(G_1)\ootimes{\opens(G_0)}\opens(G_1)\ar[rrrr]^-{m_!}&&&&\opens(G_1)\;,
}
\]
as the following diagram chase shows:
\[
\xymatrix{
(a\otimes b)\otimes c\ar@{|->}[d]_{\cong}\ar@{|->}[rrrr]^{m_!\otimes\ident}&&&&(ab)\otimes c\ar@{|->}[dd]^{m_!}\\
a\otimes (b\otimes c)\ar@{|->}[d]_{\ident\otimes m_!}\\
a\otimes(bc)\ar@{|->}[rrr]^{m_!}&&&a(bc)\ar@{=}[r]&(ab)c\;.
}
\]
Taking the right adjoints of all the above morphisms gives us the frame version of the associativity of $m$,
\[
\xymatrix{
(G_1\ptimes{G_0}G_1)\ptimes{G_0}G_1\ar[rrrr]^{m\times\ident}\ar[d]_{\cong}&&&&G_1\ptimes{G_0}G_1\ar[dd]^{m}\\
G_1\ptimes{G_0}(G_1\ptimes{G_0}G_1)\ar[d]_{\ident\times m}\\
G_1\ptimes{G_0}G_1\ar[rrrr]^{m}&&&&G_1\;,
}
\]
because, similarly to what we have argued for $u$, we have
\begin{eqnarray*}
(m\times\ident)_!&=&m_!\otimes\ident\\
(\ident\times m)_!&=&\ident\otimes m_!\;. \qed
\end{eqnarray*}
\end{proof}

Now we see a fundamental example of multiplicative quantal frame.

\begin{theorem}
Let $S$ be an abstract complete pseudogroup. Then
$\lcc(S)$ is a multiplicative quantal frame.
\end{theorem}

\begin{proof}
The right adjoint of the quantal multiplication,
\[\mu_*:\lcc(S)\to \lcc(S)\otimes_{\spp \lcc(S)}\lcc(S)\;,\]
is given by the formula
\[\mu_*(U)=\V\{V\otimes W \st VW\subseteq U\}\;.\]
Due to the universal property of $\lcc(S)$ as a sup-lattice, the question of whether $\lcc(S)$ is multiplicative, that is of whether $\mu_*$ preserves joins, is equivalent to asking whether the map
\[f:S\to\lcc(S)\otimes_{\spp \lcc(S)}\lcc(S)\]
defined by
\[f(x)=\V\{\downsegment y\otimes \downsegment z\st yz\le x\}\]
preserves all the joins that exist, in which case $\mu_*$ is the unique homomorphic extension of $f$; that is, for each compatible set $X\subseteq S$ we need to see that
\[f(\V X)\subseteq \V f(X)\;.\]
Equivalently, we need to see that $yz\le\V X$ implies $\downsegment y\otimes \downsegment z\subseteq \V f(X)$ for all $y,z\in S$.

Consider then $y,z\in S$ such that $yz\le\V X$.
For each $x\in X$ we have $yzx^{-1}x\le x$ because
\[yzx^{-1} x\le(\V X)x^{-1}x=\V_{w\in X} w x^{-1} x\le x\;,\]
where the equality is a consequence of distributivity, and the last inequality follows from the fact that $X$ is compatible and therefore
$wx^{-1}\in E(S)$ for all $w\in X$. Hence,
\[y(y^{-1}yzx^{-1}x)=(yy^{-1}y)zx^{-1}x=yzx^{-1}x\le x\]
and thus, by definition of $f$, we obtain
\[\downsegment y\otimes \downsegment (y^{-1}yzx^{-1}x)\subseteq f(x)\;.\]
From here, using distributivity, it follows that
\begin{eqnarray*}\downsegment y\otimes \downsegment \left({y^{-1}yz\V_{x\in X} x^{-1}x}\right) &=&
\downsegment y\otimes \downsegment \left({\V_{x\in X} y^{-1}yzx^{-1}x}\right)\\
&=& \downsegment y\otimes \left({\V_{x\in X} \downsegment (y^{-1}yzx^{-1}x)}\right)\\
&=& \V_{x\in X}\downsegment y\otimes \downsegment(y^{-1}yzx^{-1}x)\\
&\subseteq& \V f(X)\;.
\end{eqnarray*}
Since $\V_{x} x^{-1}x=(\V X)^{-1}(\V X)$ [this is stated in \cite[p.\ 27, Prop.\ 17]{Lawson} for $X\neq\emptyset$, but the proof applies to any $X$], we further conclude that
\[yz\V_{x\in X}x^{-1} x=yz\left({\V X}\right)^{-1}\left({\V X}\right)=yz\]
(the last equality follows from the fact that for any elements $a$ and $b$ of an inverse semigroup the condition $a\le b$ implies $a=ab^{-1}b$),
and thus
\begin{eqnarray*}
\downsegment y\otimes \downsegment z&=&\downsegment (yy^{-1} y)\otimes \downsegment z=\downsegment y\downsegment(y^{-1} y)\otimes \downsegment z=\downsegment y\otimes \downsegment(y^{-1}y)\downsegment z\\
&=&\downsegment y\otimes \downsegment(y^{-1}yz)=\downsegment y\otimes\downsegment\left({y^{-1}yz\V_{x\in X}x^{-1} x}\right)\subseteq \V f(X)\;.\qed
\end{eqnarray*}
\end{proof}

\subsection{Inverse quantal frames}\label{sec:invqufr}

Now we shall prove some facts about those quantal frames that are also inverse quantales. In particular we shall see that such quantal frames are necessarily multiplicative and of the form $\lcc(S)$, up to isomorphism, and that this gives us a category which is equivalent to the category of abstract complete pseudogroups $\apgcat$.

\begin{definition}
By an \emph{inverse quantal frame} $Q$ will be meant a supported quantal frame whose maximum is a join of partial units:
\[1=\V \ipi(Q)\;.\]
\end{definition}

We remark that any inverse quantal frame $Q$ is an inverse quantale in the sense of our original definition, due to distributivity: if $a\in Q$ then
\[a=a\wedge 1=a\wedge\V\ipi(Q)=\V \{a\wedge s\st s\in\ipi(Q)\}\;,\]
where each $a\wedge s$ is of course a partial unit. Hence, in particular, any inverse quantal frame is a stable quantal frame.

Recall the adjunction $\lcc\dashv\ipi$ of \ref{cor:adjunction}, between the category $\apgcat$ of abstract complete pseudogroups and the category $\ssq$ of stably supported quantales. For each stably supported quantale $Q$ we shall denote by
\[\varepsilon_Q:\lcc(\ipi(Q))\to Q\]
the corresponding component of the co-unit of the adjunction. This is a homomorphism of unital involutive quantales that is explicitly defined by
\[\varepsilon_Q(U)=\V U\;.\]

\begin{lemma}
Let $Q$ be an inverse quantal frame. Then $\varepsilon_Q$ is a surjective frame homomorphism whose restriction to the set of principal ideals of $\lcc(\ipi(Q))$ is injective.
\end{lemma}

\begin{proof}
$\varepsilon_Q$ is surjective because $Q$ is an inverse quantale. Hence, it remains to show that $\varepsilon_Q$ preserves binary meets. First, $\ipi(Q)$ is an abstract complete pseudogroup and thus in particular it is a meet semilattice. The fact that $\varepsilon_Q$ preserves binary meets is now an essentially immediate consequence of the coverage theorem for frames \cite{Johnstone} (with minor adaptations due to the possible absence of a maximum in $\ipi(Q)$), but a direct proof using the explicit formula for the co-unit is also immediate and we give it here: for any $U,V\in\lcc(\ipi(Q))$ we have
\[U\cap V=\{s\wedge t\st s\in U,\ t\in V\}\;,\]
and thus using the frame distributivity of $Q$ we obtain
\begin{eqnarray*}
\varepsilon_Q(U\cap V)&=&\varepsilon_Q(\{s\wedge t\st s\in U,\ t\in V\})=\V_{s\in U,\ t\in V} s\wedge t\\
&=&\V U\wedge\V V=\varepsilon_Q(U)\wedge\varepsilon_Q(V)\;.
\end{eqnarray*}
Finally, the restriction of $\varepsilon_Q$ to principal ideals is the assignment
\[\downsegment s\mapsto \V\downsegment s=s\;,\]
and, of course, this is an order embedding. \qed
\end{proof}

\begin{lemma}\label{lem:principalideals}
Let $S$ be an abstract complete pseudogroup. The principal ideals of $\lcc(S)$ form a downwards closed set.
\end{lemma}

\begin{proof}
Let $s\in S$, and let $U\in\lcc(S)$ be such that $U\subseteq\downsegment s$. For all $t,u\in U$ we have 
$tu^{-1}\le ss^{-1}\le e$, and, similarly, $t^{-1}u\le e$; that is, $t$ and $u$ are compatible and we conclude that $U$ is a compatible subset of $S$. Since $S$ is complete the join $\V U$ exists in $S$, and since $U$ is closed under joins it must contain $\V U$. Hence, $U$ is the principal ideal $\downsegment \V U$. \qed
\end{proof}

\begin{theorem}\label{thm:invqufr}
Let $Q$ be an inverse quantal frame. Then there is an isomorphism
\[\lcc(\ipi(Q))\cong Q\]
of unital involutive quantales.
\end{theorem}

\begin{proof}
$\varepsilon_Q$ is a homomorphism of unital involutive quantales, and the fact that it is an isomorphism follows immediately from \ref{prop:framebasis} and the previous two lemmas. \qed
\end{proof}

\begin{corollary}
Any inverse quantal frame is multiplicative.
\end{corollary}

Denoting by $\iqf$ the full subcategory of $\uiq$ whose objects are the inverse quantal frames, we have:

\begin{theorem}\label{thm:equivapgiqf}
The categories $\apgcat$ and $\iqf$ are equivalent.
\end{theorem}

\begin{proof}
From \ref{thm:invqufr} it follows that the adjunction $\lcc\dashv\ipi$ between $\apgcat$ and $\iqcat$ restricts to a reflection between $\apgcat$ and $\iqf$. Now let $S$ be an abstract complete pseudogroup. The unit of the adjunction gives us the injective homomorphism of abstract complete pseudogroups
\[\eta_S:S\to\ipi(\lcc(S))\]
defined by $s\mapsto\downsegment s$, and in order to prove that the reflection is in fact an equivalence it remains to see that $\eta_S$ is surjective. Let then $U\in\ipi(\lcc(S))$. By definition of partial  unit this is an element of $\lcc(S)$ such that $UU^*\subseteq E(S)$ and $U^*U\subseteq E(S)$. Hence, $st^{-1}\in E(S)$ and
$s^{-1}t\in E(S)$ for all $s,t\in U$, which means that $U$ is a compatible subset of $S$. Hence, again as in \ref{lem:principalideals}, $U$ must coincide with the principal ideal $\downsegment \V U$; that is, $\eta_S(\V U)=U$, showing that $\eta_S$ is surjective. \qed
\end{proof}

\subsection{Groupoids from quantales}\label{sec:groupoidsfromquantales}

Now we determine the conditions under which the category associated to a multiplicative quantal frame $Q$ is a groupoid. As we shall see, this happens if and only if $Q$ is an inverse quantal frame. A first step is given by the following result, which basically produces a straightforward translation of the inversion law of groupoids into the language of quantales.

\begin{lemma}\label{lemma:groupoidiffinversion}
Let $Q$ be a multiplicative quantal frame. The localic category $(G_1,G_0,d,r,u,m,i)$ associated to $Q$ is a groupoid (with inversion $i$) if and only if $Q$ satisfies the following two conditions, for all $a\in Q$:
\begin{eqnarray}
(a\wedge e)1&=&\V_{xy^*\le a} x\wedge y\;,\label{frameinversionlawleft}\\
1(a\wedge e)&=&\V_{x^*y\le a} x\wedge y\;.\label{frameinversionlawright}
\end{eqnarray}
\end{lemma}

\begin{proof}
The groupoid inversion law is the commutativity of the following diagram:
\[
\xymatrix{G_1\ar[d]_d\ar[rr]^-{\langle\ident,i\rangle}&&G_1\ptimes{G_0}G_1\ar[d]^m&&G_1\ar[ll]_-{\langle i,\ident\rangle}\ar[d]^r\\
G_0\ar[rr]_u&&G_1&&G_0.\ar[ll]^u}
\]
Consider its dual frame version:
\begin{equation}\label{frameinversionlaw}
\vcenter{\xymatrix{\opens(G_1)\ar@{<-}[d]_{d^*}\ar@{<-}[rr]^-{\lbrack\ident,i^*\rbrack}&&\opens(G_1\ptimes{G_0}G_1)\ar@{<-}[d]^{m^*}&&\opens(G_1)\ar@{<-}[ll]_-{\lbrack i^*,\ident\rbrack}\ar@{<-}[d]^{r^*}\\
\opens(G_0)\ar@{<-}[rr]_{u^*}&&\opens(G_1)&&\opens(G_0)\ar@{<-}[ll]^{u^*}}}
\end{equation}
The commutativity of the left square of (\ref{frameinversionlaw}) is equivalent, for each $a\in Q=\opens(G_1)$, to the equation
\begin{equation}\label{frameinversionlaw2}
d^*(u^*(a))=\lbrack\ident,i^*\rbrack(m^*(a))\;.
\end{equation}
Taking into account the following formulas, for all $b,x,y\in Q$,
\begin{eqnarray*}
i^*(b)&=&b^*\\
d^*(u^*(b))&=&(b\wedge e)1\\
m^*(b)&=&\V\{x\otimes y\st xy\le b\}\\
\lbrack f,g\rbrack(x\otimes y)&=&f(x)\wedge g(y)\;,
\end{eqnarray*}
we see that (\ref{frameinversionlaw2}) is equivalent to
\[(a\wedge e)1=\V_{xy\le a} x\wedge y^*\;,\]
which is equivalent to (\ref{frameinversionlawleft}).
Similarly, the right square of the diagram (\ref{frameinversionlaw}) is equivalent to (\ref{frameinversionlawright}). \qed
\end{proof}

The following two lemmas are motivated by the equations (\ref{frameinversionlawleft}) and (\ref{frameinversionlawright}). We remark that, even though the equations have been introduced in the context of multiplicative quantal frames, the lemmas hold for more general quantal frames.

\begin{lemma}
Let $Q$ be a stable quantal frame. Then for all $a\in Q$ the following inequalities hold:
\begin{eqnarray}
(a\wedge e)1&\ge&\V_{xy^*\le a} x\wedge y\label{lemma:inverseequivin1ge} \\
1(a\wedge e)&\ge&\V_{x^*y\le a} x\wedge y\label{lemma:inverseequivin2ge}\;.
\end{eqnarray}
\end{lemma}

\begin{proof}
Recall the property (\ref{prop:s10a}) of supported quantales:
\[\spp(x\wedge y)\le xy^*\;.\]
The support $\spp$ coincides with the sup-lattice homomorphism
\[u_!\circ d_!:Q\to Q\;,\]
and thus from (\ref{prop:s10a}) we obtain, by adjointness, $x\wedge y\le d^*(u^*(xy^*))$. This is equivalent to the statement that $x\wedge y\le d^*(u^*(a))$ for all $a\in Q$ such that $xy^*\le a$, and thus we obtain
\[d^*(u^*(a))\ge\bigvee_{xy^*\le a}x\wedge y\;.\]
Then (\ref{lemma:inverseequivin1ge}) is a consequence of this and of the equality $(a\wedge e)1=d^*(u^*(a))$, and (\ref{lemma:inverseequivin2ge}) is proved analogously taking into account that
$1(a\wedge e)=r^*(u^*(a))$ and, again using (\ref{prop:s10a}),
$u_!(r_!(x\wedge y))=u_!(d_!((x\wedge y)^*))=\spp(x^*\wedge y^*)\le x^*y$. \qed
\end{proof}

This has as a consequence that the conditions (\ref{frameinversionlawleft}) and (\ref{frameinversionlawright})
are equivalent, for any stable quantal frame, to the following lax version of them,
\begin{eqnarray}
(a\wedge e)1&\le&\V_{xy^*\le a} x\wedge y\label{lemma:inverseequivin1le} \\
1(a\wedge e)&\le&\V_{x^*y\le a} x\wedge y\label{lemma:inverseequivin2le}\;,
\end{eqnarray}
leading us to the second lemma:

\begin{lemma}\label{inversionlemma}
Let $Q$ be a unital involutive quantal frame. Then $Q$ satisfies the two conditions (\ref{lemma:inverseequivin1le}) and (\ref{lemma:inverseequivin2le}) if and only if $\V\ipi(Q)=1$.
\end{lemma}

\begin{proof}
Let us assume that $\V\ipi(Q)=1$ and prove from there that $Q$ satisfies (\ref{lemma:inverseequivin1le}):
\begin{eqnarray*}
(a\wedge e)1&=&(a\wedge e)\V\{x\st xx^*\le e\textrm{ and }x^*x\le e\}\\
&\le&(a\wedge e)\V\{x\st xx^*\le e\}\\
&=&\V\{(a\wedge e)x\st xx^*\le e\}\\
&\le&\V\{(a\wedge e)x\st xx^*a\le a\}\\
&=&\V\{(a\wedge e)^* x\st xx^*a\le a\}\\
&\le&\V\{a^*x\wedge x\st xx^*a\le a\}\\
&=&\V\{a^*x\wedge x\st x(a^* x)^*\le a\}\\
&\le&\V\{x\wedge y\st xy^*\le a\}\;.
\end{eqnarray*}
Proving (\ref{lemma:inverseequivin2le}) is done in an analogous way, but in the beginning retaining the inequality $x^*x\le e$ instead of $xx^*\le e$.

For the converse let us assume that both (\ref{lemma:inverseequivin1le}) and (\ref{lemma:inverseequivin2le}) hold, and from there let us prove that $\V\ipi(Q)=1$. From (\ref{lemma:inverseequivin1le}) we obtain
\begin{eqnarray*}
1=(e\wedge e)1&\le&\V\{x\wedge y\st xy^*\le e\}\\
&\le&\V\{x\wedge y\st (x\wedge y)(x\wedge y)^*\le e\}=\V\{x\st xx^*\le e\}\;,
\end{eqnarray*}
and, similarly, from (\ref{lemma:inverseequivin2le}) we obtain
\[1\le\V\{y\st y^*y\le e\}\;.\]
Hence,
\begin{eqnarray*}
1&\le&\V\{x\st xx^*\le e\}\wedge\V\{y\st y^*y\le e\}\\
&=&\V\{x\wedge y\st xx^*\le e\textrm{ and }y^*y\le e\}\\
&\le&\V\{x\wedge y\st (x\wedge y)(x\wedge y)^*\le e\textrm{ and }(x\wedge y)^*(x\wedge y)\le e\}\\
&=&\V\{x\st xx^*\le e\textrm{ and }x^*x\le e\}=\V\ipi(Q)\;. \qed
\end{eqnarray*}
\end{proof}

We finally arrive at the main result of this section.

\begin{theorem}\label{inversiontheorem}
The following conditions are equivalent.
\begin{enumerate}
\item $Q$ is an inverse quantal frame.
\item $Q$ is a multiplicative quantal frame and the category associated to $Q$ is a groupoid.
\item $Q$ is a multiplicative quantal frame and the category associated to $Q$ is an \'{e}tale groupoid.
\end{enumerate}
\end{theorem}

\begin{proof}
We already know that inverse quantal frames are multiplicative, and it is clear from the previous three lemmas that the category associated to a multiplicative quantal frame $Q$ is a groupoid if and only if $Q$ is an inverse quantal frame. What remains to be proved is therefore that this groupoid is necessarily \'{e}tale. Let then $Q$ be an inverse quantal frame. Then $\ipi(Q)$ is a cover, and, since we already know that $d$ is open, in order to show that $d$ is a local homeomorphism it suffices to prove,
for each $a\in\ipi(Q)$, that the frame homomorphism
\[f=((-)\wedge a)\circ d^*:\spp Q\to\downsegment a\]
is surjective (\cf\ section \ref{sec:locqu}). Let $x\in\downsegment a$. Then $xx^*\in\spp Q$ because $x$ is a partial unit. The conclusion that $f$ is surjective follows from the fact that $f(xx^*)=x$:
\begin{itemize}
\item By (\ref{prop:s6a}) we have $xx^*1=x1$. Then $x\le xx^*1$, and thus $x\le xx^*1\wedge a=f(xx^*)$;
\item By \ref{lemma:defstab}-\ref{prop:s15a} we have $f(xx^*)=xx^*1\wedge a=xx^*a\le xa^*a\le xe=x$. \qed
\end{itemize}
\end{proof}

\begin{definition}
Let $Q$ be an inverse quantal frame. We denote its associated localic \'{e}tale groupoid by $\groupoid(Q)$.
\end{definition}

\subsection{Separating examples}\label{sec:examples}

We have studied various kinds of quantales, and in particular we have obtained the following inclusions:
\begin{equation}\label{hierarchy}
\begin{array}{ccccc}
\left\{\begin{minipage}{1.8cm}{\center 
inverse\\
quantales\\
}
 \end{minipage}\right\}
 &&\subset&&
 \left\{\begin{minipage}{1.9cm}{\center 
stably\\
supported\\
quantales\\
}
 \end{minipage}\right\}\\
 \cup&&&&\cup\\
\left\{\begin{minipage}{1.5cm}{\center 
inverse\\
quantal\\
frames\\
}
 \end{minipage}\right\}&\subset&
\left\{\begin{minipage}{2.6cm}{\center 
multiplicative\\
quantal\\
frames\\
}
 \end{minipage}\right\}&\subset&
\left\{\begin{minipage}{1.5cm}{\center 
stable\\
quantal\\
frames\\
}
 \end{minipage}\right\}\;.
 \end{array}
\end{equation}
The examples that follow show that all the inclusions are strict.

\begin{example}\label{fuzzyexample}
A stable quantal frame which is not multiplicative is the commutative quantale $Q=\pwset X$ with
\begin{itemize}
\item $X=\{1,x\}$ (with $1\neq x$),
\item trivial involution,
\item $\spp U=\{1\}$ for all $U\neq\emptyset$,
\item $e=\{1\}$,
\item multiplication defined on the atom $\{x\}$ by the condition $\{x\}\{x\}=\{1,x\}$ (and freely extended to unions in each variable).
\end{itemize}
This also shows that not every stably supported quantale is an inverse quantale.
\end{example}

\begin{example}
Let $M$ be the idempotent ordered monoid which, besides the unit $1$, contains only one additional element $x$ such that $1\le x$. The set $Q=\lc(M)$ of downwards closed subsets of $M$ is an idempotent unital quantale under pointwise multiplication, with $e=\{1\}$, and it is commutative. With trivial involution, and with a support defined by $\spp(U)=\{1\}$ for all $U\neq\emptyset$, we obtain a multiplicative quantal frame that is not an inverse quantale because $M$ is not a union of partial units.

We remark that this example differs from the previous one because the quantale is not a powerset.
There is a good reason for this: if the powerset of (the set of arrows of) a small discrete involutive category is a supported quantale at all, then it is necessarily an inverse quantale. This is because the axiom $\spp a\le aa^*$ alone forces the category to be a groupoid, since on singletons the axiom gives us $d(x)=\spp(\{x\})\subseteq \{xx^*\}$, \ie, $d(x)=xx^*$, and also $r(x)=d(x^*)=x^*x$, and thus the involution of the category coincides with inversion: $x^*=x^{-1}$.
\end{example}

\begin{example}
Now we present an example of an inverse quantale that is not a frame, hence showing that the two vertical inclusions of the diagram (\ref{hierarchy}) are strict.
Consider a non-T$_0$ topological space $X=\{x,y,z\}$ with three open sets $\emptyset$, $\{x,y\}$, and $X$. This space has exactly one non-idempotent automorphism, namely the bijection $s$ that permutes $x$ and $y$. Hence, the order structure of its pseudogroup $S=\ipi(X)$ is as follows:
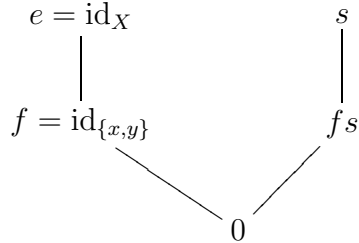
\[
\xymatrix{
e=\ident_X\ar@{-}[d]&&s\ar@{-}[d]\\
f=\ident_{\{x,y\}}\ar@{-}[dr]&&fs\ar@{-}[dl]\\
&0
}
\]
The multiplication of $S$ is commutative, it is defined by the conditions $s^2=e$ and $E(S)=\{0,f,e\}$, and each element is its own inverse.
The inverse quantal frame $Q=\lcc(S)$, which in this example coincides with $\lc(S)$, has nine elements (we write $s$ instead of $\downsegment s$, $f\vee s$ instead of $\downsegment f\cup\downsegment s$, etc.), namely $0$, $f$, $t=fs$, $e$, $a=f\vee fs$, $s$, $b=e\vee fs$, $c=f\vee s$, and $1=e\vee s$, where of course we have $a=f1$. It is now straightforward to obtain the multiplication table of $Q$; we present only the upper triangle because $Q$ is commutative:
\[
\begin{array}{ccccccccccc}
&\vline&0&f&t&e&a&s&b&c&1\\
\cline{0-10} 0&\vline&0&0&0&0&0&0&0&0&0\\
f&\vline&&f&t&f&a&t&a&a&a\\
t&\vline&&&f&t&a&f&a&a&a\\
e&\vline&&&&e&a&s&b&c&1\\
a&\vline&&&&&a&a&a&a&a\\
s&\vline&&&&&&e&c&b&1\\
b&\vline&&&&&&&b&c&1\\
c&\vline&&&&&&&&b&1\\
1&\vline&&&&&&&&&1
\end{array}
\]
Now consider the equivalence relation $\theta$ on $Q$ whose only non-singular equivalence class is $\{b,c,1\}$. The rightmost three entries of each line of the table are always equivalent, which means that $\theta$ is a congruence for the multiplication. Similarly, any join of an element of $Q$ with either $b$, $c$ or $1$ necessarily produces an element in $\{b,c,1\}$ (because $b$ and $c$ are maximal elements of $Q$), and thus $\theta$ is also a congruence for binary joins (and hence for all joins because $Q$ is finite). Since $\theta$ is trivially also a congruence for the involution, the quotient $Q/{\theta}$ is a unital involutive quantale with seven elements ordered as follows:
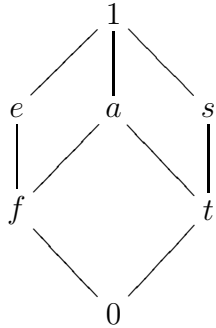
\[
\xymatrix{
&1\ar@{-}[dl]\ar@{-}[dr]\ar@{-}[d]\\
e\ar@{-}[d]&a\ar@{-}[dl]\ar@{-}[dr]&s\ar@{-}[d]\\
f\ar@{-}[dr]&&t\ar@{-}[dl]\\
&0
}
\]
This lattice is not distributive (for instance we have $s\wedge(e\vee a)=s\wedge 1=s$ and $(s\wedge e)\vee (s\wedge a)=0\vee t=t\neq s$), but it is a supported quantale because $\theta$ is a congruence also with respect to the support, since in $Q$ we have $\spp b=\spp c=\spp 1=e$. Hence, $Q/\theta$ is an inverse quantale but not an inverse quantal frame.
\end{example}

\section{Groupoid quantales}\label{section:quantalgroupoids}\label{sec:groupoids}

In this section we describe some applications of quantal frame techniques to localic and topological groupoids. In particular, we shall obtain new characterizations of \'{e}tale groupoids. The results of this section can also be seen as a partial converse to those of the previous one, since we are now concerned with studying quantal frames that are obtained from groupoids, and they establish an equivalence between the concepts of localic \'{e}tale groupoid and inverse quantal frame.

\subsection{Quantal groupoids}

As we have mentioned in section \ref{sec:introgrpdqnts}, if the topology $\topology(G)$ of a topological groupoid $G$ is closed under pointwise multiplication of open sets then $\topology(G)$ is a unital involutive quantale.
The localic analogue of this is of course a localic groupoid $G$
whose multiplication map $m$ is open, but even just by assuming that $m$ is semiopen relevant conclusions are obtained. In particular, as we shall see below, in that case $\opens(G_1)$ is a quantal frame, which motivates the following definition.

\begin{definition}
By a \emph{quantal groupoid} is meant a localic groupoid whose multiplication map is semiopen. If $G$ is a quantal groupoid, the \emph{groupoid quantale} of $G$, denoted by $\opens(G)$, is defined to be the involutive quantale of the following theorem.
\end{definition}

\begin{theorem}
Let $G$ be a localic groupoid:
\[
\xymatrix{
G_1\ptimes{G_0}G_1\ar[r]^-m&G_1\ar@(ru,lu)[]_i\ar@<1.2ex>[rr]^r\ar@<-1.2ex>[rr]_d&&G_0.\ar[ll]|u
}
\]
If $G$ is quantal then $\opens (G_1)$ is a quantale whose multiplication
\[m':\opens(G_1)\otimes\opens(G_1)\to\opens(G_1)\]
is the sup-lattice homomorphism $m'=m_!\circ q$, where
\[q:\opens(G_1\times G_1)\to\opens(G_1\ptimes{G_0} G_1)\]
is the frame quotient that defines $G_1\times_{G_0} G_1$ as a sublocale of $G_1\times G_1$. Furthermore, this quantale has an involution given by
\[a^*=i_!(a)=i^*(a)\;.\]
\end{theorem}

\begin{proof}
The proof of associativity of $m'$ is, with direction reversed, entirely analogous to the proof of associativity in \ref{theorem:localiccategory}.

Let us prove that $i_!$ is an involution on $\opens(G)$. The first condition, namely $a^{**}=a$, follows from $i\circ i=\ident$ (\cf\ \ref{prop:invcat}), as does the fact that $i_!=i^*$. For the second condition, $(ab)^*=b^* a^*$, we begin by recalling the equation
\[
i\circ m=m\circ\chi\;,
\]
where $\chi$ is the isomorphism $\langle i\circ\pi_2,i\circ\pi_1\rangle$, which satisfies $\chi_!=\chi^*$ because $\chi\circ\chi=\ident$ (\cf\ \ref{prop:invcat}). In particular, we have
\begin{eqnarray*}
\chi_!(a\otimes b)&=&\chi^*(a\otimes b)=\lbrack \pi_2^*\circ i^*,\pi_2^*\circ i^*\rbrack(a\otimes b)
=\pi_2^*(i^*(a))\wedge\pi_1^*(i^*(b))\\
&=&1\otimes a^*\wedge b^*\otimes 1=b^*\otimes a^*\;.
\end{eqnarray*}
Hence, noting that $i_!\circ m_!=m_!\circ\chi_!$, we obtain
\[(ab)^*=i_!(m_!(a\otimes b))=m_!(\chi_!(a\otimes b))=m_!(b^*\otimes a^*)=b^* a^*\;.\qed \]
\end{proof}

\begin{lemma}
Let $G$ be a quantal groupoid. Then $d^*(a)$ is right-sided and $r^*(a)$ is left-sided, for all $a\in\opens(G_0)$.
\end{lemma}

\begin{proof}
Let $a\in\opens(G_0)$. We have, in $\opens(G)$,
\begin{equation}\label{eq:qugrpd}
d^*(a)1=m_!(d^*(a)\otimes 1)=m_!(\pi_1^*(d^*(a)))\;,
\end{equation}
where $\pi_1:G_1\times_{G_0}G_1\to G_1$ is the first projection. One of the defining conditions of $G$ as a localic groupoid is $d\circ\pi_1=d\circ m$, and thus we can replace $\pi_1^*$ by $m^*$ in (\ref{eq:qugrpd}), which leads to
\[d^*(a)1=m_!(m^*(d^*(a)))\le d^*(a)\;.\]
In a similar way one proves that $r^*(a)$ is left-sided. \qed
\end{proof}

\begin{lemma}\label{lemma:rightadjointtosupport}
Let $G$ be a quantal groupoid. Then, for all $a\in\opens(G)$, we have
\begin{eqnarray}
d^*(u^*(a))&=&\V_{bc^*\le a}b\wedge c\label{lemma:rightadjointtosupport1}\\
r^*(u^*(a))&=&\V_{b^*c\le a}b\wedge c\label{lemma:rightadjointtosupport2}\;.
\end{eqnarray}
\end{lemma}

\begin{proof}
This is entirely similar to the proof of \ref{lemma:groupoidiffinversion}, where the groupoid inversion law was seen to be equivalent to the two conditions (\ref{frameinversionlawleft}) and (\ref{frameinversionlawright}), except that now we cannot assume equations like $d^*(u^*(a))=(a\wedge e)1$ or $r^*(u^*(a))=1(a\wedge e)$, which make no sense because we do not even have a unit $e$. \qed
\end{proof}

\begin{lemma}\label{lemma:semiopendr}
Any quantal groupoid has semiopen domain and range maps, with $d_!(a)=u^*(a1)$ and $r_!(a)=u^*(1a)$.
\end{lemma}

\begin{proof}
Let us verify the conditions $d_! d^*\le \ident$ and $d^* d_!\ge\ident$ with respect to the proposed definition of $d_!$ (for $r$ it is analogous). Let $a\in\opens(G_0)$. Taking into account that $d^*(a)$ is right-sided we obtain
\[d_!(d^*(a))=u^*(d^*(a)1)\le u^*(d^*(a))=a\;.\]
Now let
$a\in\opens(G_1)$. Then
$d^*(d_!(a))$ equals $d^*(u^*(a1))$ which, by (\ref{lemma:rightadjointtosupport1}), equals
\[\V\{b\wedge c\st bc^*\le a1\}\;,\]
and this is greater or equal to $a$ (for instance, let $b=c=a$). \qed
\end{proof}

\begin{example}\label{exm:opengroupoids}
Examples of quantal groupoids are common. For instance, observing that the multiplication $m$ of any groupoid can be obtained as a pullback of the domain map $d$ along itself,
\[
\vcenter{\xymatrix{X\ar@/_/[ddr]_g\ar@/^/[drrr]^f\ar@{.>}[dr]|h\\
&G_1\ptimes{G_0}G_1\ar[rr]^{\pi_1}\ar[d]^m&&G_1\ar[d]^d\\
&G_1\ar[rr]_d&&G_0}}\hspace*{1cm}h=\langle f, m\circ\langle i\circ f,g\rangle\rangle
\]
one concludes that any open groupoid, by which we mean a localic groupoid with $d$ open, is an example of a quantal groupoid, and one whose multiplication is open rather than just semiopen, because open maps are stable under pullback.

The converse also holds, \ie, assuming that $m$ is open we conclude that $d$ is open (by \ref{lemma:semiopendr} we know that it is semiopen) because it satisfies the Frobenius reciprocity condition:
\[
\begin{array}{rcll}
d_!(a\wedge d^*(b))&=&u^*((a\wedge d^*(b))1)&\textrm{(By \ref{lemma:semiopendr}.)}\\
&=& u^*(m_!((a\wedge  d^*(b))\otimes 1))\\
&=& u^*(m_!(\pi_1^*(a\wedge d^*(b)))) &\textrm{($\pi_1^*(-)=-\otimes 1$\,.)}\\
&=& u^*(m_!(\pi_1^*(a)\wedge \pi_1^* d^*(b)))\\
&=& u^*(m_!(\pi_1^*(a)\wedge m^* d^*(b)))&\textrm{($d\circ m=d\circ \pi_1$\,.)}\\
&=& u^*(m_!(\pi_1^*(a))\wedge d^*(b))&\textrm{(Frobenius cond.\ for $m$\,.)}\\
&=& u^*(a1\wedge d^*(b)) &\textrm{($\pi_1^*(a)=a\otimes 1$\,.)}\\
&=& u^*(a1)\wedge u^* d^*(b)\\
&=& d_!(a)\wedge b&\textrm{(By \ref{lemma:semiopendr} and $d\circ u=\ident$\,.)}
\end{array}
\]
\end{example}

\subsection{Localic \'{e}tale groupoids}\label{sec:locetfgrps}

\begin{definition}
A localic groupoid $G$ is said to be \emph{unital} if the map \[u:G_0\to G_1\] is open (and thus $G_0$ is an open sublocale of $G_1$).
\end{definition}

\begin{lemma}
Let $G$ be a unital quantal groupoid. The involutive quantale $\opens(G)$ is unital, and the multiplicative unit is $e=u_!(1_{G_0})$.
\end{lemma}

\begin{proof}
The proof is the same, with direction reversed, as the proof of the unit laws in \ref{theorem:localiccategory}. \qed
\end{proof}

\begin{lemma}
Let $G$ be a unital quantal groupoid. The following conditions hold in $\opens(G)$, for all $a\in\opens(G_0)$ and $b\in\opens(G_1)$:
\begin{eqnarray}
b\wedge d^* (a)&=& u_! (a) b \label{bda}\\
b\wedge r^*(a) &=& b u_!(a) \label{bra}
\end{eqnarray}
\end{lemma}

\begin{proof}
Let $a\in\opens(G_0)$. We have
\[e\wedge r^*(a)= u_! (1_{G_0})\wedge r^*(a)=u_!(1_{G_0}\wedge u^*( r^*( a)))=u_! (u^* (r^* (a)))=u_! (a)\;,\]
where the second equality follows from the Frobenius reciprocity condition for $u$, and the last equality is a consequence of the condition $r\circ u=\textrm{id}$. Using this, and recalling from the end of section \ref{sec:toplocgrpds} that the pushout of $d^*$ and $r^*$
\[\opens(G_1)\ootimes{\opens(G_0)}\opens(G_1)\]
satisfies the condition
\begin{equation}\label{pushoutcondition}
b\otimes (d^*(a)\wedge c) = (b\wedge r^*(a))\otimes c\;,
\end{equation}
we prove (\ref{bda}):
\[b\wedge d^* (a)=m_!(e\otimes(b\wedge d^*(a)))=m_!((e\wedge r^*(a))\otimes b)=m_!(u_!(a)\otimes b)=u_!(a)b\;.\]
Equation (\ref{bra}) is proved in a similar way, this time starting from the condition $e\wedge d^*(a)=u_!(a)$, which is an instance of (\ref{bda}). \qed
\end{proof}

\begin{lemma}\label{lemma:ud}

Let $G$ be a unital quantal groupoid. Then, for all $a\in\opens(G_1)$,
\begin{eqnarray*}
d^*(u^*(a))&=&(a\wedge e)1\\
r^*(u^*(a))&=&1(a\wedge e)\;.
\end{eqnarray*}
\end{lemma}

\begin{proof}
Let $a\in\opens(G_1)$. From (\ref{bda}) we have
\[d^*(u^*(a))=1\wedge d^*(u^*(a))=u_!(u^*(a))1\;.\]
And we have $u_!(u^*(a))=a\wedge e$ because $u$ is open:
\[u_!(u^*(a))=u_!(u^*(a)\wedge 1_{G_0})=a\wedge u_!(1_{G_0})=a\wedge e\;.\]
Hence, $d^*(u^*(a))=(a\wedge e)1$.
For $r^*(u^*(a))$ we use (\ref{bra}) and everything is analogous. \qed
\end{proof}

\begin{theorem}\label{thm:quantalesfromgroupoids}
Let $G$ be a unital quantal groupoid. Then $\opens(G)$ is an inverse quantal frame, and its groupoid $\groupoid(Q)$ is isomorphic to $G$.
\end{theorem}

\begin{proof}
First, we show that the sup-lattice homomorphism
\[\spp = u_!\circ d_!\]
defines a support:
\begin{itemize}
\item $\spp a=u_!(d_!(a))\le u_!(1_{G_0})=e$, which proves (\ref{def:s1}).
\item An instance of (\ref{lemma:rightadjointtosupport1}) gives us
\[ d^*(u^*(aa^*))=\V_{xy^*\le aa^*}x\wedge y\;,\]
and thus $a\le d^*(u^*(aa^*))$ (make $a=x=y$). Hence, by adjointness we obtain
$u_!(d_!(a))\le aa^*$, \ie, we have proved (\ref{def:s2}).
\item Now we prove (\ref{def:s3}):
\[\begin{array}{rcll}
\spp a a&=& m_!(\spp a\otimes a)\\
&=& m_!(u_! (d_!(a))\otimes a)\\
&=& m_!((e\wedge r^*(d_!(a)))\otimes a) &\textrm{[By (\ref{bra})]}\\
&=& m_!(e\otimes(d^*(d_!(a))\wedge a)) &\textrm{[By (\ref{pushoutcondition})]}\\
&=& d^*(d_!(a))\wedge a\\
&=& a &\textrm{($d^*\circ d_!\ge \ident$)}\;.
\end{array}
\]
\end{itemize}
Hence, $\opens(G)$ is a supported quantal frame. Furthermore, from \ref{lemma:rightadjointtosupport} and \ref{lemma:ud} we obtain the equations
\begin{eqnarray*}
(a\wedge e)1 &=& \V_{xy^*\le a} x\wedge y\\
1(a\wedge e) &=& \V_{x^* y\le a} x\wedge y\;,
\end{eqnarray*}
which show, by \ref{inversionlemma}, that $\opens(G)$ is an inverse quantal frame. Hence, $\opens(G)$ has an associated (\'{e}tale) groupoid $\groupoid(\opens(G))$. Let us denote this by $\widehat{G}$, with structure maps $\hat d$, $\hat r$, $\hat u$, $\hat m$, and $\hat i$. We shall prove that $G$ and $\widehat G$ are isomorphic. Since obviously we have $G_1=\widehat{G}_1$, it is natural to look for an isomorphism $(f_1,f_0):G\to \widehat G$ with $f_1=\ident_{G_1}$. Then $f_0$ must be given by $f_0=\hat d\circ f_1\circ u=\hat d\circ u$, and, similarly, its inverse must be given by $d\circ \hat u$.
Let us verify that the pair $(f_1,f_0)$ commutes with $d$ and $\hat d$, \ie, that the following diagram commutes:
\begin{equation}\label{f0f1diagramford}
\xymatrix{
G_1\ar[d]_d\ar[rrr]^{f_1=\ident}&&&\widehat{G}_1=G_1\ar[d]^{\hat d}\\
G_0\ar[rrr]_{f_0=\hat{d}\circ u}&&&\widehat{G}_0
}
\end{equation}
In order to do this, first we remark that the results of section \ref{sec:sqf} on stable quantal frames give for the homomorphisms ${\hat d}^*\circ{\hat u}^*$ and ${\hat r}^*\circ{\hat u}^*$ the formulas
\begin{eqnarray*}
({\hat d}^*\circ {\hat u}^*)(a)&=&(a\wedge e)1\\
({\hat r}^*\circ {\hat u}^*)(a)&=&1(a\wedge e)\;,
\end{eqnarray*}
which are identical to those of \ref{lemma:ud} for $d^*\circ u^*$ and $r^*\circ u^*$, thus yielding the following identities of locale maps:
\begin{eqnarray}
\hat u\circ \hat d&=&u\circ d \label{eq:udud}\\
\hat u\circ \hat r&=&u\circ r\;.\label{eq:urur}
\end{eqnarray}
Hence, using (\ref{eq:udud}) we have
\[f_0\circ d=\hat d\circ u\circ d=\hat d\circ \hat u\circ \hat d=\ident\circ \hat d=\hat d\circ\ident=\hat d\circ f_1\;;\]
that is, the diagram (\ref{f0f1diagramford}) commutes.
Using again (\ref{eq:udud}) we show that $(f_1,f_0)$ commutes with $u$ and $\hat u$,
\[\hat u\circ f_0=\hat u\circ \hat d\circ u=u\circ d\circ u=u\circ\ident=\ident\circ u=f_1\circ u\;,\]
and using (\ref{eq:urur}) we conclude that $(f_1,f_0)$ commutes with $r$ and $\hat r$:
\[f_0\circ r=\hat d\circ u\circ r=\hat d\circ \hat u\circ \hat r=\ident\circ \hat r=\hat r\circ\ident=\hat r\circ f_1\;.\]
Hence, $(f_1,f_0)$ is a morphism of reflexive graphs. The fact that it is an isomorphism with
$f^{-1}_0= d\circ \hat u$ follows again from (\ref{eq:udud}):
\[
\begin{array}{c}(\hat d\circ u)\circ (d\circ \hat u)=\hat d\circ \hat u\circ \hat d\circ \hat u=\ident\circ\ident=\ident\;,\\
(d\circ \hat u)\circ (\hat d\circ u)=d\circ u\circ d\circ u=\ident\circ\ident=\ident\;.
\end{array}\]
From these results it follows that the pullback of $d$ and $r$ coincides with the pullback of $\hat d$ and $\hat r$ (both pullbacks are, as frame pushouts, given by the same quotient of $\opens(G_1)\otimes\opens(G_1)$), and thus it is obvious that $m=\hat m$, since both $m^*$ and ${\hat m}^*$ are right adjoint to the same quantal multiplication $\opens(G_1\times_{G_0}G_1)\to \opens(G_1)$. Similarly, $i=\hat i$ because both $i^*$ and ${\hat i}^*$ coincide with the quantale involution, and we conclude that $(f_1,f_0)$ is an isomorphism of groupoids. \qed
\end{proof}

Our results provide equivalent but new alternative definitions for the notion of \'{e}tale groupoid:

\begin{corollary}\label{cor:quantalunitalisetale}
For any localic groupoid $G$, the following are equivalent:
\begin{enumerate}
\item $G$ is \'{e}tale.
\item $G$ is quantal and unital.
\item $G$ is open and unital.
\end{enumerate}
\end{corollary}

\begin{proof}
$2\Rightarrow 1$: Immediate consequence of \ref{inversiontheorem} and \ref{thm:quantalesfromgroupoids}.\\
$3\Rightarrow 2$: Immediate because being open implies being quantal (\cf\ \ref{exm:opengroupoids}).\\
$1\Rightarrow 3$: For an \'{e}tale groupoid all the structure maps are local homeomorphisms, and thus, in particular, both $m$ and $u$ are open. \qed
\end{proof}

We remark that as consequences of \ref{inversiontheorem} and
\ref{thm:quantalesfromgroupoids} we have obtained a duality between \'{e}tale groupoids and inverse quantal frames, which is given by isomorphisms
\begin{eqnarray*}
G&\cong&\groupoid(\opens(G))\\
Q&\cong&\opens(\groupoid(Q))\;.
\end{eqnarray*}
[Indeed there is an obvious equality $Q=\opens(\groupoid(Q))$.]
It is natural to ask how well this duality behaves with respect to morphisms, a question that we shall briefly address now.

\begin{lemma}\label{lemma: morphisms}
Let $Q_1$ and $Q_2$ be inverse quantal frames, and $f:Q_1\to Q_2$ a sup-latttice homomorphism. Let also $m_1$ and $m_2$ be the multiplication maps of $\groupoid(Q_1)$ and $\groupoid(Q_2)$, respectively. The following conditions are equivalent:
\begin{enumerate}
\item $f(a)f(b)\le f(ab)$ for all $a,b\in Q_1$.
\item $(f\otimes f)\circ m_1^*\le m_2^*\circ f$.
\end{enumerate}
\end{lemma}

\begin{proof}
$1\Rightarrow 2$: Condition $1$ is equivalent to $(m_2)_!\circ (f\otimes f)\le f\circ (m_1)_!$, which by adjointness is equivalent to $f\otimes f\le m_2^*\circ f\circ (m_1)_!$. Composing with $m_1^*$ on the right we obtain $(f\otimes f)\circ m_1^*\le m_2^*\circ f\circ (m_1)_!\circ m_1^*$, and thus
we obtain condition $2$ because $(m_1)_!\circ m_1^*\le \ident$.\\
$2\Rightarrow 1$: From condition $2$ we obtain, by adjointness, $(m_2)_!\circ(f\otimes f)\circ m_1^*\le f$. Now composing, on both sides of this inequality, with $(m_1)_!$ on the right we obtain condition $1$ because $m_1^*\circ (m_1)_!\ge\ident$. \qed
\end{proof}

\begin{theorem}\label{thm:assignment}
Let $G$ and $G'$ be \'{e}tale groupoids, and let \[h=(h_1,h_0):G'\to G\] be a morphism of groupoids. Then for all $a,b\in\opens(G)$ we have \[h_1^*(a) h_1^*(b)\le h_1^*(ab)\;.\]
\end{theorem}

\begin{proof}
This a corollary of the previous lemma, since a morphism $(h_1,h_0)$ of groupoids preserves multiplication and that is equivalent to the equality $(h_1^*\otimes h_1^*)\circ m^*=(m')^*\circ h_1^*$. \qed
\end{proof}

The following example shows that the inequality in the above theorem is in general not an equality.

\begin{example}\label{exm:laxhom}
As an example of an \'{e}tale groupoid consider a nontrivial discrete group $G$ (written multiplicatively), and let $h:G\to G$ be an endomorphism. In general $h^{-1}:\pwset G\to\pwset G$ is not a homomorphism of quantales: for instance, if $h(g)=1$ for all $g\in G$, and $U=V^{-1}=\{g\}$ with $g\neq 1$, then $h^{-1}(U)=h^{-1}(V)=\emptyset$, whence
$h^{-1}(U)h^{-1}(V)=\emptyset$ but $h^{-1}(UV)=h^{-1}(\{1\})=G$.
\end{example}

This shows that there is no immediate contravariant functor from \'{e}tale groupoids to quantales, and that in order to find a duality between \'{e}tale groupoids and inverse quantal frames in the categorical sense one must be willing to change the morphisms under consideration, for instance allowing more general homomorphisms of quantales, in particular homomorphisms that are lax on multiplication as in \ref{thm:assignment}, or restricting consideration of maps of groupoids to those whose inverse images preserve quantale multiplication, etc. (Similar problems apply to multiplicative units --- for instance, in \ref{exm:laxhom} we have $e=\{1\}$ and $h^*(e)=\ker h$, and thus $h^*(e)=e$ if and only if $h$ is injective.)

\subsection{Topological groupoids}

We shall now obtain for topological groupoids some results that are analogous to those just obtained for localic groupoids. In fact some of them are in a sense more general because even if the topology of a topological groupoid is not a quantale under pointwise multiplication, the open sets can nevertheless be multiplied, and formulas similar to those of localic groupoids arise (something analogous may exist for a localic groupoid $G$ if one defines the product of two sublocales $L$ and $M$ of $G_1$ to be the image in $G_1$ of
$L\times_{G_0} M\to G_1\times_{G_0}G_1\stackrel m\to G_1$). The following lemma, whose proof is analogous to that of \ref{lemma:groupoidiffinversion}, provides the first example of this, where formulas that formally coincide with (\ref{lemma:inverseequivin1le}) and (\ref{lemma:inverseequivin2le}) are obtained. The proof techniques used in the present section are borrowed from previous parts of this paper, but the results obtained are completely independent. To simplify we shall adopt, for any topological groupoid $G$, the conventions $G=G_1$ and $G_0\subseteq G$.

\begin{lemma}
Let $G$ be any topological groupoid. Then we have, for any open set $U\in\topology(G)$:
\begin{eqnarray}
(U\cap G_0)G&\subseteq&\bigcup\{X\cap Y\st X,Y\in\topology(G),\ XY^{-1}\subseteq U \}
\label{lemma:inverseequivin1letop}\\
G(U\cap G_0)&\subseteq&\bigcup\{X\cap Y\st X,Y\in\topology(G),\ X^{-1}Y\subseteq U \}
\label{lemma:inverseequivin2letop}
\end{eqnarray}
\end{lemma}

\begin{proof}
Recall the groupoid inversion law,
\[
\xymatrix{G\ar[d]_d\ar[rr]^-{\langle\ident,i\rangle}&&G\ptimes{G_0}G\ar[d]^m&&G\ar[ll]_-{\langle i,\ident\rangle}\ar[d]^r\\
G_0\ar[rr]_u^{\subseteq}&&G&&G_0,\ar[ll]^u_{\supseteq}}
\]
which on inverse images gives us the following commutative diagram:
\begin{equation}\label{inviminversionlaw}
\vcenter{\xymatrix{\topology(G)\ar@{<-}[d]_{d^{-1}}\ar@{<-}[rr]^-{\langle\ident,i\rangle^{-1}}&&\topology(G\ptimes{G_0}G)\ar@{<-}[d]^{m^{-1}}&&\topology(G)\ar@{<-}[ll]_-{\langle i,\ident\rangle^{-1}}\ar@{<-}[d]^{r^{-1}}\\
\topology(G_0)\ar@{<-}[rr]_{u^{-1}}&&\topology(G)&&\topology(G_0)\ar@{<-}[ll]^{u^{-1}}}}
\end{equation}
The commutativity of the left square of (\ref{inviminversionlaw}) is equivalent, for each $U\in \topology(G)$, to the equation
\begin{equation}\label{inviminversionlaw2}
d^{-1}(u^{-1}(U))=\langle\ident,i\rangle^{-1}(m^{-1}(U))\;.
\end{equation}
Now consider the following obvious formulas, for all $V,X,Y\in \topology(G)$:
\begin{eqnarray*}
i^{-1}(V)&=&V^{-1}\\
d^{-1}(u^{-1}(V))&=&(V\cap G_0)G\\
m^{-1}(V)&=&\{(x, y)\in G\times_{G_0}G\st xy\in V\}\\
\langle\ident,i\rangle^{-1}(X\times_{G_0}Y)&=&\{x\in G\st x\in X,\ x^{-1}\in Y\}\\
&=&X\cap Y^{-1}\;.
\end{eqnarray*}
Since $G\times_{G_0}G$ has a basis of open sets of the form $X\times_{G_0}Y$ with $X,Y\in\topology(G)$, from the continuity of $m$ it follows that for any $(x,y)\in G\times_{G_0}G$ and any open neighborhood $V$ of $xy=m(x,y)$ there is a basic open $X\times_{G_0} Y$ containing $(x,y)$ such that $m(X\times_{G_0} Y)\subseteq V$; that is, such that $XY\subseteq V$. Hence, the open set
\[m^{-1}(V)=\{(x, y)\in G\times_{G_0}G\st xy\in V\}\]
can be rewritten as
\[m^{-1}(V)=\bigcup\{X\times_{G_0}Y\st X,Y\in \topology(G),\ XY\subseteq V\}\;,\]
and thus from (\ref{inviminversionlaw2}) we obtain (\ref{lemma:inverseequivin1letop}). Similarly, the right square of the diagram (\ref{inviminversionlaw}) gives us (\ref{lemma:inverseequivin2letop}). \qed
\end{proof}

From here, applying \ref{inversionlemma}, it immediately follows that if the topology $\topology(G)$ of a topological groupoid $G$ is a unital quantale under pointwise multiplication, with unit $G_0$, then $\ipi(G)$ is an open cover of $G$. But, in fact, even just by assuming that $G_0$ is open we arrive at the same conclusion. The proof is analogous to that of (the relevant half of) \ref{inversionlemma}:

\begin{lemma}\label{lemma:rdiscrete}
Let $G$ be any topological groupoid whose unit subspace $G_0\subseteq G_1=G$ is open. Then $\ipi(G)$, the set of open G-sets of $G$, is a cover of $G$.
\end{lemma}

\begin{proof}
Using the hypothesis that $G_0$ is open, from (\ref{lemma:inverseequivin1letop}) we obtain:
\begin{eqnarray*}
G=(G_0\cap G_0)G&\subseteq&\bigcup\{X\cap Y\st X,Y\in\topology(G),\ XY^{-1}\subseteq G_0\}\\
&\subseteq&\bigcup\{X\cap Y\st X,Y\in\topology(G),\ (X\cap Y)(X\cap Y)^{-1}\subseteq G_0\}\\
&=&\bigcup\{X\in\topology(G)\st XX^{-1}\subseteq G_0\}\;.
\end{eqnarray*}
Similarly, from (\ref{lemma:inverseequivin2letop}) we obtain
\[G\subseteq\bigcup\{Y\in\topology(G)\st Y^{-1}Y\subseteq G_0\}\;.\]
Hence,
\begin{eqnarray*}
G&\subseteq&\bigcup\{X\in\topology(G)\st XX^{-1}\subseteq G_0\}\cap\bigcup\{Y\in\topology(G)\st Y^{-1}Y\subseteq G_0\}\\
&=&\bigcup\{X\cap Y\st X,Y\in\topology(G),\ XX^{-1}\subseteq G_0,\ Y^{-1}Y\subseteq G_0\}\\
&\subseteq&\bigcup\{X\in\topology(G)\st XX^{-1}\subseteq G_0\textrm{ and }X^{-1}X\subseteq G_0\}=\bigcup\ipi(G)\;. \qed
\end{eqnarray*}
\end{proof}

\begin{theorem}\label{thm:topetale}
Let $G$ be a topological groupoid. The following conditions are equivalent:
\begin{enumerate}
\item $G_0$ is open and $d$ is an open map.
\item $G$ is \'{e}tale.
\item $G_0$ is open and the topology $\topology(G)$ is closed under pointwise multiplication of open sets.
\item $G_0$ is open and the pointwise product $UG$ is open for each $U\in\topology(G)$.
\item $G_0$ is open and the pointwise product $UU^{-1}$ is open for each $U\in\topology(G)$.
\end{enumerate}
\end{theorem}

\begin{proof}
($1\Rightarrow 2$) If $G_0$ is open then, by \ref{lemma:rdiscrete}, $\ipi(G)$ is an open cover of $G$, and thus for each $x\in G$ there is an open neighborhood $U\in\ipi(G)$ of $x$ such that the restriction of $d$ to $U$ is injective. Since $d$ is also open, it follows that its restriction to $U$ is a homeomorphism from $U$ onto the open set $d(U)$. Hence, $d$ is a local homeomorphism.

($2\Rightarrow 3$) If a groupoid is \'{e}tale then both its multiplication map $m$ and the unit inclusion map $u$ are local homeomorphisms, hence open.

($3\Rightarrow 4$ and $3\Rightarrow 5$) Trivial.

($4\Rightarrow 1$ and $5\Rightarrow 1$) Under either of the hypotheses 4 or 5 the domain map is open because its direct image satisfies the following equations, for all subsets $U\subseteq G$ (these, of course, are the formulas for the stable support of $\pwset G$):
\[d(U)=UG\cap G_0=UU^{-1}\cap G_0\;. \qed\]
\end{proof}

\subsection{Inverse semigroups as \'{e}tale groupoids}\label{sec:sasg}

One of the consequences of the results so far is that we have been provided, as remarked in section 1, with an entirely new way in which to construct localic and topological \'{e}tale groupoids: what one really constructs is an inverse quantal frame, the groupoid being a derived object. Furthermore, since the categories involved are not equivalent, this new kind of construction is different in a non trivial way from constructions that are carried out ``within'' the category of groupoids. We shall now give a very simple example of this, one which will also provide insights into the nature of the groupoids that can be associated to a fundamental example of inverse quantal frame, namely $\lc(S)$ for an inverse semigroup $S$.

The most immediate way of obtaining a (discrete) groupoid from an inverse semigroup $S$ is to consider the idempotents to be the units, and, for $x,y\in E(S)$, the arrows $s:x\to y$ to be those $s\in S$ such that $ss^{-1}=x$ and $s^{-1} s=y$. Let us denote this groupoid by $\widetilde S$. The passage from  $S$ to $\widetilde S$ implies that information is lost, and keeping track of this information is one of the motivations behind the notion of ordered groupoid (see \cite{Lawson}). Of course, a partial order can be equivalently described by means of an Alexandroff topology (or, in this case, a co-Alexandroff topology), and thus the following theorem can be regarded as a topological reinterpretation of the ordered groupoid of an inverse semigroup, obtained via quantales.

\begin{theorem}\label{thm:sgasgrpd}
Let $S$ be an inverse semigroup. Then $\widetilde S$, equipped with $\lc(S)$ as a topology, is an \'{e}tale groupoid, and the inverse quantale structure that this groupoid induces on $\lc(S)$ is the same as in \ref{lemma:invqfrominvs}.
\end{theorem}

\begin{proof}
The domain map of $\widetilde S$ is defined by $d(s)=ss^{-1}$. This is monotone with respect to the natural order of $S$, and thus continuous with respect to the topologies $\topology(\widetilde S)$ and $\topology(\widetilde S_0)$ (the latter being the subspace topology $\lc(E(S))$).
Similarly, the inversion map $s\mapsto s^{-1}$ is continuous. The space of composable pairs $\widetilde S\times_{\widetilde S_0}\widetilde S$ is again equipped with the co-Alexandroff topology (with respect to the direct product order), and thus the multiplication map of $\widetilde S$ is continuous because it is monotone. Hence, $\widetilde S$ is a topological groupoid.

Now let $U,V\in\topology(\widetilde S)$. Let $UV$ denote their multiplication in the sense of \ref{lemma:invqfrominvs}, and let $U\& V$ denote the pointwise multiplication induced by the groupoid multiplication of $\widetilde S$. It is clear that we have $U\& V\subseteq UV$. Let us check that also the converse holds. Let $x\in UV$. Then $x=yz$ for some $y\in U$ and $z\in V$, where the product $yz$ is given by the semigroup multiplication. But $yz=yy^{-1}yzz^{-1}z=yfz$, where $f=y^{-1}yzz^{-1}$ is idempotent, and thus $yz$ equals the product $(yf)(fz)$. Now we have
$yf\le y$ and $fz\le z$, and thus $yf\in U$ and $fz\in V$. Also, $r(yf)=(yf)^{-1}(yf)=fy^{-1}yf=fy^{-1}y=y^{-1}yzz^{-1}y^{-1}y=y^{-1}yzz^{-1}=f$ and, similarly, $d(fz)=f$, and thus the groupoid multiplication $U\& V$ contains the product $(yf)(fz)$, which is $x$. Hence, $U\& V=UV$. Besides showing that the inverse quantale structure of $\lc(S)$ coincides with that of $\topology(\widetilde S)$, this allows us to apply \ref{thm:topetale}, since in addition $E(S)$ is an open set, and conclude that $\widetilde S$ is an \'{e}tale groupoid. \qed
\end{proof}

\begin{corollary}
The inverse semigroup $S$ can be recovered from the topological groupoid $\widetilde S$. The multiplication of $S$ is the multiplication in $\lc(S)$ of the principal ideals.
\end{corollary}

Another aspect that this example brings out is that the construction of $\widetilde S$ can be understood in terms of a universal property, as already mentioned in section \ref{sec:introduction}: the natural topology of this groupoid is $\lc(S)$ and it is ``freely'' generated by $S$. Hence, $\widetilde S$ is in a suitable sense a ``universal groupoid'' of $S$. Furthermore, $\lc(S)$ is an inverse quantal frame and thus $\groupoid(\lc(S))$ is a localic \'{e}tale groupoid whose spectrum is, in an obvious sense, the soberification of the topological groupoid $\widetilde S$. The nature of this spectrum, including the relation it has to the notion of germ of an element of an inverse semigroup, depends again crucially on the properties of $\lcc(S)$ and will be addressed elsewhere.

\end{document}